\documentclass[sn-mathphys-num,iicol]{sn-jnl}

\usepackage{graphicx}%
\usepackage{multirow}%
\usepackage{amsmath,amssymb,amsfonts}%
\usepackage{amsthm}%
\usepackage{mathrsfs}%
\usepackage[title]{appendix}%
\usepackage{xcolor}%
\usepackage{textcomp}%
\usepackage{manyfoot}%
\usepackage{booktabs}%
\usepackage{algorithm}%
\usepackage{algorithmicx}%
\usepackage{algpseudocode}%
\usepackage{listings}%
\usepackage{float}

\theoremstyle{thmstyleone}%
\theoremstyle{thmstyletwo}%

\theoremstyle{thmstylethree}%

\usepackage{graphicx} 
\usepackage{algorithm}  
\usepackage{algpseudocode}
\usepackage{color}
\usepackage{setspace}
\usepackage{subcaption}
\usepackage{siunitx}

\usepackage{listings}
\usepackage[section]{placeins}
\usepackage{hyperref}

\usepackage{etoolbox}
\makeatletter
\patchcmd{\ttlh@hang}{\parindent\z@}{\parindent\z@\leavevmode}{}{}
\patchcmd{\ttlh@hang}{\noindent}{}{}{}

\makeatother
\usepackage{amsfonts}
\usepackage{bm}
\usepackage{bbm}
\usepackage{chngcntr}
\usepackage{indentfirst}
\usepackage{graphicx}

\usepackage{mwe}

\usepackage[export]{adjustbox}
\usepackage{graphicx}
\usepackage{caption}
\usepackage{subcaption}

\usepackage{booktabs,caption}
\usepackage[font=small,skip=5pt]{caption}
\usepackage{multirow}
\usepackage[section]{placeins}

\def\vc#1{\mathbf{\boldsymbol{#1}}}     
\def\tn#1{\boldsymbol{#1}}

\def \E{{\mathbb E}}

\def\avg#1{\langle#1\rangle}
\def \D{{{\rm I\kern-.3em D}}}

\def\avg#1{\langle#1\rangle}

\def\vc#1{\mathbf{\boldsymbol{#1}}}     
\def\tn#1{\boldsymbol{#1}}

\def\avg#1{\langle#1\rangle}

\def\std{\mathop{\rm std}}

\def\div{\operatorname{div}}
\def\grad{\nabla}

\def\vmu{\vc\mu}

\begin{document}

\title[Article Title]{Deep learning surrogate for predicting hydraulic conductivity tensors from DFM models}


\author*[1]{\fnm{Martin} \sur{Špetlík}}\email{martin.spetlik@tul.cz}

\author[1]{\fnm{Jan} \sur{Březina}}\email{jan.brezina@tul.cz}

\author[2]{\fnm{Eric} \sur{Laloy}}\email{eric.laloy@sckcen.be}

\affil*[1]{\orgdiv{Institute of New Technologies and Applied Informatics, Faculty of Mechatronics, Informatics and
Interdisciplinary Studies}, \orgname{Technical University of Liberec}, \orgaddress{\street{Studentská 1402/2}, \city{Liberec}, \postcode{461 17}, \country{Czech Republic}}}

\affil[2]{\orgdiv{Engineered and Geosystems Analysis Unit, Sustainable Waste Management and Decommissioning}, \orgname{Belgian Nuclear Research Center (SCK CEN)}, \orgaddress{\street{Boeretang 200}, \city{Mol}, \postcode{2400}, \country{Belgium}}}



\abstract{
Simulating water flow in fractured crystalline rock requires tackling its stochastic nature. We aim to utilize the multilevel Monte Carlo method for cost-effective estimation of simulation statistics. This multiscale approach entails upscaling of fracture hydraulic conductivity by homogenization.
In this work, we replace 2D numerical homogenization based on the discrete fracture-matrix (DFM) approach with a surrogate model to expedite computations. We employ a deep convolutional neural network (CNN) connected to a deep feed-forward neural network as the surrogate. The equivalent hydraulic conductivity tensor $\tn K^{eq}$ is predicted based on the input tensorial spatial random fields (SRFs) of hydraulic conductivities, along with the cross section and hydraulic conductivity of fractures.
Three independent surrogates with the same architecture are trained, each with a different ratio of fracture to matrix hydraulic conductivity $K_f/K_m$.
As the ratio $K_f/K_m$ increases, the distribution of $\tn K^{eq}$ becomes more complex, leading to a decline in the prediction accuracy of the surrogates.
The prediction accuracy improves as the fracture density decreases, regardless of the $K_f/K_m$. We also investigate prediction accuracy for different correlation lengths of input SRFs.
The observed speedup gained by surrogates varies from $4\times$ to $28\times$ depending on the number of homogenization blocks.
Upscaling by numerical homogenization and surrogate modeling is compared on two macroscale problems. For the first one, the accuracy of outcomes is directly correlated with the precision of $\tn K^{eq}$ predictions. For the latter one, we observe only a mild impact of the upscaling method on the accuracy of the results.}




\maketitle

\section{Introduction}
Modeling groundwater processes is integral to researching deep geological repositories of radioactive waste. The repository planned in the Czech Republic will be located in crystalline rock, where a small intergranular pore space causes groundwater to flow predominantly through a fracture network \cite{Banks2002}.
Given the inherent uncertainties in the modeled processes, our aim is to utilize the multilevel Monte Carlo method (MLMC) \cite{Giles2015} to obtain statistics at a reduced computational cost compared to the plain Monte Carlo method.
%
Applying MLMC directly to models with discrete fractures, such as discrete fracture networks (DFN), poses challenges. Coarser-level models in MLMC cannot capture fractures below their resolution, and DFNs themselves struggle to resolve fractures smaller than about 1/1000th of the domain size due to an exponential growth in fracture count. The impact of the small fractures, however, should not be ignored. To address this limitation, we adopt a discrete fracture-matrix (DFM) approach that integrates both continuum and discrete fractures. This allows us to homogenize the impact of small fractures into the matrix properties when transitioning from a finer to a coarser scale DFM model.
%
A surrogate based on deep learning techniques can accelerate a costly numerical homogenization procedure. 
This study investigates deep learning-based surrogates for the numerical homogenization of an equivalent hydraulic conductivity tensor for 2D DFM models.

Homogenization techniques typically utilize a representative
elementary volume (REV) to calculate the effective hydraulic conductivity of fractured media \cite{Auriault20090101}. However, REV is not applicable for fracture networks in rocks due to their fractal structure 
\cite{Bonnet2001ScalingOF}. The homogenized hydraulic conductivity becomes a heterogeneous random field with an increased lower bound for correlation length. In such cases, the terminology shifts to equivalent properties instead of effective properties (see \cite{https://doi.org/10.1029/2005RG000169, RENARD1997253}).
For our work, we utilize numerical homogenization in a block-averaged sense, where the correlation length of the resulting field varies depending on the block size (see \cite{CHEN201560}). While some authors differentiate between homogenization and upscaling, e.g., \cite{https://doi.org/10.1002/fld.267}, we use them interchangeably in this work.


Our objective is to apply numerical homogenization to calculate an equivalent hydraulic conductivity tensor.
The determination of the equivalent hydraulic conductivity of fractured porous media has been reviewed in \cite{geosciences12070269, RENARD1997253, https://doi.org/10.1029/2005RG000169}. Additionally, numerical homogenization in the context of DFM models has also been studied. 
Bogdanov et al. \cite{https://doi.org/10.1029/2001WR000756} examine a cubic rock matrix containing single-size two-dimensional polygonal fractures. The rock matrix is assumed to be isotropic and homogeneous. As an extension of the work by Koudina et al. \cite{Koudina1998}, they employ two distinct permeabilities: one for the fractures and another for the rock matrix. 
A tetrahedral mesh with fractures is utilized to model the Darcy flow separately for the porous matrix and the fracture network. Bogdanov and colleagues also investigate the impact of fracture density and matrix hydraulic conductivity on the upscaled hydraulic conductivity. 
Their subsequent work \cite{PhysRevE.76.036309} extends the previous work by considering power-law distributed fracture sizes.
Lang et al. \cite{Lang20140815} determine an equivalent hydraulic conductivity tensor using a pressure gradient and a flux averaging technique for both 2D and 3D DFM models on unstructured grids. Other relevant studies include, for example, \cite{https://doi.org/10.1029/2000WR900340, Azizmohammadi2017}.

Deep learning (DL) \cite{Goodfellow-et-al-2016} is an extensive part of machine learning that utilizes deep architectures of artificial neural networks.
In our case, deep neural networks (DNNs) are employed to solve regression problems using data-driven supervised learning~\cite{Goodfellow-et-al-2016}. An unknown mapping of input data to output data is learned from collected sample pairs of inputs and corresponding outputs of a DFM model. 
The input comprises a spatial random field (SRF) on an unstructured mesh, along with a discrete fracture network with random properties. The output corresponds to an equivalent hydraulic conductivity tensor.
Deep learning methods have recently proven successful in determining the properties of porous media. Specifically, supervised convolutional neural networks (CNNs) have gained popularity for this purpose.
Studies such as \cite{CAGLAR2022106973}, \cite{ALQAHTANI2020106514} or \cite{WANG2022123284} adopt 2D CNNs to predict permeability or porosity based on input images of a porous media. 
Many researchers \cite{RAO2020109850,  Hong2020, Peng2022PhNetPM, VASILYEVA2021185, https://doi.org/10.1029/2022JB025378} focus on more realistic 3D models, where 3D CNNs are employed.
Hong and Liu \cite{Hong2020} train four 3D CNNs for four different target permeabilities. They consider 3D binary images of the pore geometry as input and four permeabilities (three axial directions and mean permeability) as the corresponding output. Meng et al. \cite{MENG2023104520} introduced a CNN with a transformer neural network to predict permeability.
To predict equivalent permeability tensors, various approaches have been employed, including deep feed-forward neural networks \cite{STEPANOV2023114980}, \cite{pr11020601}, graph convolutional neural networks \cite{Cai_2023}, and generative adversarial networks \cite{FERREIRA2022104264}.
Although He et al. \cite{10.2118/203901-MS} and Andrianov \cite{Andrianov2022UpscalingOT} use CNNs to upscale DFM models, they consider constant hydraulic conductivity for both the matrix and fractures. This approach is inadequate for our needs, as we intend to incorporate our surrogate into stochastic calculations using MLMC.

Our work represents a pioneering effort in employing surrogates to predict equivalent hydraulic conductivity tensors from 2D DFM models, where both matrix and fracture elements exhibit non-constant hydraulic conductivity. Our study on 2D surrogates serves as a proof of concept before dealing with the more applicable but more complex 3D homogenization case.
Our contribution involves the utilization of tensorial spatial random fields (SRFs) for matrix conductivities and non-constant hydraulic conductivity of fractures, which depends on fracture size and aperture. 
We investigate the impact of the matrix-to-fracture hydraulic conductivity ratio on the distribution of the equivalent hydraulic conductivity tensor and the associated predictive capabilities of the trained surrogates. 
The accuracy of the trained surrogates is further analyzed using unseen DFM model data with varying SRF correlation lengths or number of fractures. 
Finally, we compare the trained surrogates to numerical homogenization using the resulting macroscale hydraulic conductivity field within two macroscale problems: groundwater flow through a given surface (Aquifer problem) and upscaling the hydraulic conductivity tensor (Anisotropy problem).

The article is structured as follows. Section \ref{dfm_model_section}
explains multiscale DFM models as a conjunction of a discrete fracture network (DFN) and an equivalent continuous media (ECM) approach. The numerical homogenization of the equivalent hydraulic conductivity tensor is described, including two investigated macroscale problems. The numerical homogenization is put into the context of multiscale models.
Section \ref{meta_model_arch_section} provides the surrogate architecture and metrics for assessing surrogate accuracy. Section \ref{dataset_config} describes the dataset configurations and preprocessing adopted. Section \ref{results_section} presents and discusses the results, including the remaining challenges.
Finally, Section \ref{conclussion_section} provides the conclusions of the article.

\def\rmin{\underline r}
\def\rmax{\overline r}
\section{Multiscale models of fractured media}\label{dfm_model_section}
 The primary aim of this section is to introduce homogenization problems, which will later be replaced by surrogates. The second aim is to contextualize these homogenization problems within a multiscale description of fractured media and the multilevel Monte Carlo method. 
 We start with a simple DFN stochastic model of the distribution of discrete fractures. Then, we proceed to the ECM description of the rock mass by a random permeability tensor field. We then combine these descriptions using a DFM approach. Refer to \cite{Berre2019} for various descriptions of fractured rock. We continue with homogenization and two macroscale test problems. Finally, we describe the motivation behind homogenization using fast surrogates to enable the multilevel Monte Carlo method for Darcy flow in fractured rocks.

\subsection{Discreate fracture network}
As Darcy flow in crystalline rock is primarily attributed to fractures, discrete fracture network (DFN) models are commonly employed: \cite{Long1982Porousa}, \cite{Bour1997Connectivity}, \cite{deDreuzy2012Influence}.
The characteristics of fractures in DFNs can include (1)~length, (2)~aperture, (3)~orientation, or (4)~density (as well as connectivity), all of which are regarded as random variables.
\begin{enumerate}
\item Based on evidence from nature, it is commonly assumed that fracture lengths ($f_s$) follow a power-law distribution (see \cite{Bonnet2001ScalingOF}):
\begin{equation}
    \label{pow_law}
    f_s \sim C r^{-\alpha},
    \qquad
    C = \frac{1-\alpha}{\rmax^{1-\alpha} - \rmin^{1-\alpha}},
\end{equation}
where $\alpha$ represents the characteristic exponent, while $\rmin$ and $\rmax$ denote the minimum and maximum lengths of a fracture, respectively.

\item The aperture of a fracture plays a crucial role in determining its permeability. According to Bonnet et al. \cite{Bonnet2001ScalingOF}, various distributions such as log-normal, power-law, normal, or Lévy-stable have been utilized to describe the aperture distribution. However, for simplification purposes, some studies (e.g., \cite{https://doi.org/10.1002/2015JB011879, Lang20140815}) adopt a constant aperture for all fractures. We assume a linear relation between fracture aperture and size: $\delta = af_{s}$, with $a=\num{1e-4}$.

\item  
The distribution of orientations is frequently modeled using von Mises or Fisher distributions (for further details, refer to \cite{Adler1999}). For the sake of simplicity, we consider an isotropic distribution of orientations.

\item 
Dimensionless density, as described in \cite{Adler1999}, serves to quantify the degree of connectivity of fractures. It determines the number of fracture centers (i.e., the total number of fractures) within a given domain size. In the 2D case, the dimensionless density is defined as follows:
\begin{equation}
\rho^{\prime}_{2D} = \rho_{2D}A_{ex},
\end{equation}
where $\rho_{2D}$ is the number of fractures per unit area, and $A_{ex}$ denotes the excluded area, i.e., an area around a fracture where the center of another fracture cannot be located to ensure the two fractures do not intersect, for details see \cite[ch. 3]{Sahimi20110420}.
The concept of dimensionless density can be extended naturally to the 3D case with 2D fractures.
\end{enumerate}

For simplicity, this article assumes uncorrelated fractures with positions determined by a Poisson process and independent orientations. See 
\cite[p. 31]{Sahimi20110420} and \cite{Liu2016} for more general models.
Additionally, we consider the cubic law for fracture hydraulic conductivity:
\begin{equation}
    \label{eq:k_frac}
    K_f = \frac{g \rho_w}{12\mu}\delta^{2},
\end{equation}
where $g$ is the gravitational acceleration, $\rho_w$ denotes the density of water, and $\mu$ stands for the dynamic viscosity of water.

\subsection{Spatially correlated random field on the matrix}\label{matrix_srf_section}
The DFN approach becomes computationally prohibitive when dealing with a large number of fine-scale fractures. 
The equivalent continuous media (ECM) approach, as detailed in \cite{Hadgu_comparative_2017} and \cite{Kottwitz_Investigating_2021}, offers a solution by representing rock hydraulic properties with a random permeability tensor field $\tn K(x)$. This random field can be characterized either implicitly, through the homogenization of a DFN realization (refer to Section \ref{sec:DFM}), or explicitly, by describing its point distribution and correlation structure. In particular, we consider the permeability tensor field in the form:
\begin{equation*}
  \tn K = Q^T \Lambda Q,
  \quad
  Q = \begin{pmatrix}
        Y_x & -Y_y \\ 
        Y_y & Y_x
    \end{pmatrix},
    \quad
  \Lambda = \begin{pmatrix}
        k_x & 0\\ 
        0 & k_y
    \end{pmatrix}.
\end{equation*}
We assume the isotropic rotation matrix $Q$, which is consistent with the isotropic distribution of fractures. In order to prescribe the rotation with a correlation structure, we use unit vector field $\vc Y = \vc X / |\vc X|$, where $\vc X = (X_x, X_y)$ has components that are independent Gaussian fields with zero mean, unit variance, and a specified correlation length $\lambda$. The principal permeabilities are given by the Gaussian vector 
$(k_x, k_y) = \exp \left( \vc\mu + \sqrt{\vc \Sigma}\,\vc{\tilde k}\right)$,
 where $\vc\mu$ is the mean vector, $\vc \Sigma$ is the covariance matrix, and $\vc{\tilde k}$ is a Gaussian field. The components of $\vc{\tilde k}$ are independent, follow a standard $N(0,1)$ distribution, and share a correlation length $\lambda$.
Consequently, the tensor field is defined by four independent scalar fields: $X_x$, $X_y$, $\tilde k_x$, and $\tilde k_y$, each following a $N(0,1)$ distribution with a Gaussian correlation of length $\lambda$. For sampling these fields, the GSTools library by \cite{GSTools} is utilized.

\subsection{Discrete fracture-matrix model}\label{sec:DFM}
The discrete fracture-matrix (DFM) approach, as outlined in  \cite{Sandve_efficient_2012a}, \cite{Berrone_Simulations_2013}, and \cite{Brezina2016Analysis}, integrates aspects of both DFN and ECM through a mixed mesh composed of both 2D and 1D (or 3D and 2D) finite elements. The computational domain is decomposed into the set of fractures and the matrix domain: $\Omega=\Omega_f \sqcup \Omega_m$. In $\Omega_m$, we consider Darcy flow equation:
\begin{equation}
  \div \delta_m\vc u_m = 0,\quad  \vc u_m = - \tn K_m \grad h_m, 
\end{equation}
where $\vc u_m(\vc x)$ denotes the Darcy velocity $[m/s]$ and $\tn K_m(\vc x)$ is the permeability tensor $[m^2]$. To emphasize the true scaling of the quantities with respect to the characteristic length $L$, we introduce a non-physical principal unknown $h_m = p_m / \mu$ $[1/s]$, which is associated with the pressure $p_m$ $[Pa]$ and the dynamic viscosity $\mu$ $[Pa \cdot s]$. 
Although the usage of permeability is good to justify the scaling, for the sake of practicality, we will henceforth use hydraulic conductivity out of this section, retaining the notation K. For the sake of simplicity of the homogenization formulas and consistency with the equation on $\Omega_f$, we set $\delta_m = 1$. 

The corresponding equation on $\Omega_f$ is given by:
\begin{equation}
  \div \delta_f \vc u_f = q^+ + q^-, \quad 
  \vc u_f =-k_f \grad h_f,
\end{equation}
where $\delta_f(\vc x)=a r_i$ $[m]$ is the aperture of fracture $i$ and  $k_f(\vc x)=\delta^2 / 12$ denotes the scalar fracture permeability. The sources $q^\pm$ denote outflows from the two aligned boundaries of $\Omega_m$ with normal vectors $\vc n^\pm$, coupling the fracture and matrix domain through Robin-like boundary conditions:
\begin{equation}
  -\tn K_m \grad h_m \cdot \vc n^\pm 
  = q^\pm := \frac{k_f}{2\delta}(h_m^\pm - h_f) \quad \text{on }\Omega_f.
\end{equation}
Our proprietary multiphysical software, Flow123d (\cite{flow123d}), is employed to solve DFM Darcy flow problems using mixed finite elements. 

We now introduce two benchmark problems designed to assess the performance of the studied homogenization surrogates. Importantly, the second benchmark problem serves a dual role, also acting as the microscale homogenization problem to be approximated by the surrogates. The multiscale nature of this approach is further elaborated in the subsequent section.
\\
\subsubsection{Aquifer problem}\label{aquifer_problem}
We examine fluid flow within a square domain $\Omega = (0, L)^2$, where no-flow boundary conditions are applied on the top and bottom sides of the domain. We impose Dirichlet boundary conditions on the vertical sides, setting the hydraulic head $h = H$ at $x=0$ and $h = 0$ at $x=L$. Let $Y$ $[m^2s^{-1}]$ denote the total horizontal flux through the domain. The equivalent horizontal permeability $k_x$ can then be calculated as:
\begin{equation}
k_{x} = \frac{Y}{H}.
\end{equation}
Clearly, $Y$ is linear with respect to the prescribed pressure gradient $H/L$, so $k_x$ is independent of $H$. Next, using unit analysis, we obtain the following scale equivariance: $k_x[L, L^2\tn K_m, L^2 k_f] = L^2 k_x[1, \tn K_m, k_f]$.
\\

\subsubsection{Anisotropy problem}\label{anisotropy_problem}
To obtain the full permeability tensor, we solve two problems on $\Omega = (0, L)^2$, applying the boundary pressure head $P^1=x$ and $P^2=y$, respectively. 
Mixed formulation with $RT_0$ finite elements provides velocity vectors $\vc u(e, P)$ constant on the elements. Using element permeability tensors, we can approximate 
the pressure gradients as $\nabla h(e, P) = -\tn K_e^{-1} \vc u(e, P)$. For a vector quantity $\vc v(e, P)$ on elements, we introduce a weighted average:
\begin{equation}
  \avg{\vc v}^i = \frac{\sum_{e\in\mathcal T} |e| \delta_e\vc v(e, P^i)}{\sum_{e\in\mathcal T} |e| \delta_e}.
\end{equation}
We fit the homogenized symmetric tensor $\tn K^{eq}$ by solving the least squares problem:
\begin{equation}
-\begin{bmatrix}
\avg{\nabla h}^1_x & \avg{\nabla h}^1_y& 0 \\
0 & \avg{\nabla h}^1_x & \avg{\nabla h}^1_y \\
\avg{\nabla h}^2_x & \avg{\nabla h}^2_y& 0 \\
0 & \avg{\nabla h}^2_x & \avg{\nabla h}^2_y \\
\end{bmatrix}
\begin{bmatrix}
k_{xx}^{eq} \\
k_{xy}^{eq} \\
k_{yy}^{eq} \\
\end{bmatrix}
=
\begin{bmatrix}
\avg{u}_{x}^1 \\
\avg{u}_{y}^1 \\
\avg{u}_{x}^2 \\
\avg{u}_{y}^2
\end{bmatrix}
\end{equation}
Note that $\avg{\nabla h}_x^1 \approx \avg{\nabla h}_y^2 \approx 1$ and $\avg{\nabla h}_y^1 \approx \avg{\nabla h}_x^2 \approx 0$, due to prescribed boundary conditions $P^1$ and $P^2$. Similarly to the previous case, the homogenization problem exhibits scale equivariance:
\begin{equation}
\label{eq:equivariance}
\tn K^{eq}[L, L^2 \tn K_m, L^2 k_f] = L^2 \tn  K^{eq}[1, \tn K_m, k_f]
\end{equation}

\subsection{Numerical homogenization in the context of multiscale models}\label{num_hom_multiscale}
\begin{figure*}[!htb]
  \centering
  \includegraphics[width=\linewidth]{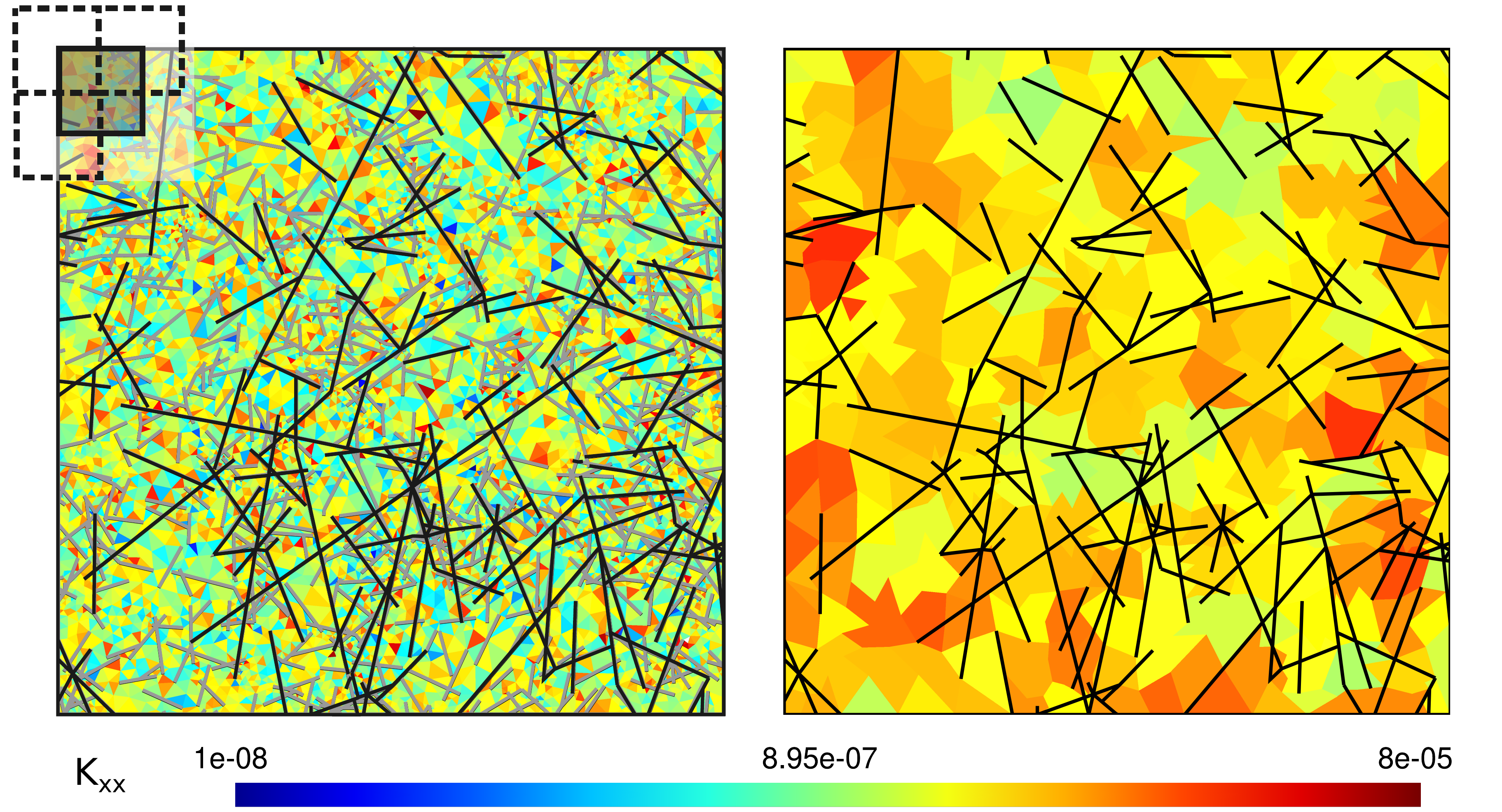}
  \centering
  \caption{{\bf Left:} The fine-scale DFM model ($h<5$), grey fractures, and permeability tensor field are homogenized using overlapping square blocks of size $15$ depicted in the corner. {\bf Right:} The coarse model ($H<10$) with homogenized permeability tensor component $K_{xx}$. Notice the reduced spread of the upscaled field.}
  \label{fig:fine_coarse_srf_dfm_model}
\end{figure*}

For the benchmark problems \ref{aquifer_problem} and~\ref{anisotropy_problem}, a fine mesh discretization can be employed to explicitly resolve all fractures within the size range $(\rmin, \rmax) = (h, L)$. A locally refined mesh with an element size
upper bound $h \ll L$ is necessary for discretization compatible with fractures, as depicted in Fig. \ref{fig:fine_coarse_srf_dfm_model}, left. 
We denote the problem result by $Q_h(\mathcal F_{h, L}, K_h)$, where $\mathcal F_{h, L}$ represents a specific realization of the fracture network for this size range, and $\tn K_h$ denotes a realization of the permeability random field effect of fractures up to size $h$.
Considering the high computational cost of fine discretizations, a two-scale approach is employed as an alternative.  This consists of a macro problem 
$Q_H(\mathcal F_{H, L}, \tn K_H)$ using a coarse mesh with maximal element size $H$, where $H > h$, and the coarse permeability tensor field $\tn K_H$ determined in terms of the micro problem \ref{anisotropy_problem}.  

In order to be compatible with convolutional neural networks, the numerical homogenization problems are solved on a grid of homogenization blocks, similar to \cite{Durlofsky1991}. Four homogenization blocks are illustrated in Fig. \ref{fig:fine_coarse_srf_dfm_model}
with center points $\vc{x}$, size $l\times l$, $l=1.5 H$, and overlap $l/2$.
In order to evaluate boundary blocks, both fracture networks and the permeability field $\tn K_h$ are generated on a slightly larger domain of side $L+2l$. 
For each block, we determine the intersection of the fractures $\mathcal F_{h, H}$ with the block, construct a compatible mesh with a maximal element size of $h$, and interpolate the field $\tn K_h$ to it. Then we calculate, the block's equivalent permeability $\tn K_H(\vc x)= P_h(\vc x, \mathcal F_{h, H}, \tn K_h)$. Finally, linear interpolation is used to transfer raster $\tn K_H$ values to the coarse grid cells, as exemplified in Fig. \ref{fig:fine_coarse_srf_dfm_model}, right.

The two-scale procedure described will serve as a testbed for simulating the application of developed neural network surrogates, specifically designed to replace the homogenization problem $P_h$.
While the acceleration of the two-scale model offered by homogenization surrogates is valuable on its own, the primary motivation lies in their application within the multilevel Monte Carlo (MLMC) method, as introduced by \cite{Giles2015}. The MLMC estimate of the expectation $\E[Q]$ is based on a telescoping sum:
\begin{equation}
\avg{Q}_{N_1,\dots, N_I} = \sum_{i=1}^I \avg{Q^{f,i} - Q^{c,i}}_{N_i},
\end{equation}
where $N_i$ denotes the number of samples at resolution level $i$ and 
$\avg{X}_{N_i}$ is the average of $X$ using $N_i$ samples. Both the fine approximation $Q^{f, i}$ and the coarse approximation $Q^{c, i}$ of the quantity $Q$ at level $i$ are random variables dependent on the same random input $\omega$. 
The direct fine discretization and the two-scale model are used as follows: 
\begin{equation}
\begin{aligned} 
Q^{f,i} &= Q_{h_i}\big(\mathcal F_{h_i, L}, \tn K_{h_i}\big),\\
Q^{c,i} &= Q_{H_i}\big(\mathcal F_{H_i, L}, 
    P_{h_i}(\mathcal F_{h_i, H_i}, \tn K_{h_i})\big),
\end{aligned}
\end{equation}
where $H_i = h_{i-1}$. All fracture sets $\mathcal F$ and fine tensor fields $\tn K_{h_i}$ are functions of a common input random state $\omega$. To obtain a functional MLMC estimator, we must ensure that:
\begin{equation}
\E Q^{f,i} = \E Q^{c, i+1}, 
\end{equation}
Discussion of this condition is beyond the scope of this paper; however, it is satisfied provided that the homogenized hydraulic conductivity $\tn K_{H_i}$ has 
the same distribution as $\tn K_{h_{i-1}}$ sampled directly. 
Regarding the requirements for the homogenization surrogates applicable to MLMC, we can make two key observations:

\begin{itemize}
    \item The homogenization problems $P_{h_i}$ only differ in the length scale determined by the block size $l_i = 1.5 H_i$. The mesh resolution and size range of discrete fractures are the same, relative to $l_i$. Consequently, we may employ scale equivariance (Eq. (\eqref{eq:equivariance}) and evaluate 
    the homogenized permeability field on the unit square: 
    $\tn K_{H_i}(x) = l_i^2 P_h(x, 1, K_{h_i}/l_i^2)$.
    \item While the fracture hydraulic conductivity $k_f$ scales with $l_i^2$ according to Eq. (\eqref{eq:k_frac}), this is not generally the case for the matrix hydraulic conductivity $\tn K_{h_i}$.  Therefore, the ratio $K_f/K_m$ varies between MLMC levels, necessitating a surrogate model that is robust to different hydraulic conductivity ratios.
\end{itemize}

\section{Surrogate architecture}\label{meta_model_arch_section}
Our deep learning (DL) surrogate solves a regression problem through supervised learning. It learns an unknown mapping from input random DFM parameters (SRF of matrix hydraulic conductivity tensors and fractures' hydraulic conductivity and cross section) to output equivalent tensors of hydraulic conductivity, utilizing collected sample pairs of inputs and corresponding outputs of a DFM model.
The DL surrogate comprises a convolutional neural network (CNN) serving as a feature extractor and a feed-forward neural network (FNN) responsible for the regression task.

Convolutional neural networks \cite[ch. 9]{Goodfellow-et-al-2016} (CNNs) process data on regular grids, particularly images. Numerous architectures of CNNs (see~\cite{Alzubaidi2021}) have recently been developed for tasks such as image classification, segmentation, natural language processing, etc.
CNNs utilize two fundamental operations: convolution and pooling. 
Popular software libraries Tensorflow and PyTorch use the term "convolution" to refer to cross-correlation.
Consider a 2D image $P$ of shape $I \times J \times C_{in}$, where $C_{in}$ stands for the number of input channels, e.g., RGB images have three channels.
Let $R$ be a convolution kernel of shape $M \times N \times C_{in} \times C_{out}$, where $C_{out}$ denotes the number of output channels. 
Then, the cross-correlation operation can be expressed as: 
\begin{multline}
S(i,j, c_{out}) = (P \ast R)(i,j, c_{out}) = \\ \sum_m \sum_n \sum_{c_{in}} P(i-m, j-n, c_{in})R(m,n, c_{in}, c_{out}). 
\end{multline}
Features extracted from the data ($S$) are passed through a non-linear activation function (e.g., ReLU), and pooling is often applied to enhance computational efficiency. A pooling layer usually reduces the dimensions of $S$ by computing summary statistics (such as max, average, or sum) over small regions of each channel of $S$ (e.g., $2 \times 2$, $4 \times 4$, etc.).
Typically, the number of learnable parameters (i.e., kernel size) varies across convolutional layers. In particular, the number of channels often increases.
Each channel can be thought of as responding to a different set of features, allowing different channels to specialize in recognizing various objects, as described by Zhang et al. \cite{zhang2020dive}.
The key advantage of CNNs lies in parameter sharing. At each layer, the convolutional kernel $R$ utilizes a single set of parameters shared across all locations. This characteristic contributes to the efficiency of CNNs compared to traditional fully connected neural networks (see \cite[p. 328]{Goodfellow-et-al-2016}, for details).

We support both unstructured and structured computational meshes and utilize the Datashader software library (see \cite{datashader}) for nearest-neighbor interpolation of potentially unstructured data into a $256\times256$ matrix. We consider tensors of hydraulic conductivity on each mesh element. For each pixel in the resulting matrix, we assign values of three independent components of the 2D tensor, along with one additional value representing the cross section of the element.

The CNN part of the surrogate consists of a Conv2D layer, which integrates convolution and max pooling. It employs a $3 \times 3$ convolutional kernel with a stride of $1$ and zero padding. The kernel is accompanied by a $2 \times 2$ max pool with a stride of $2$. We use batch normalization to facilitate the training process.
After passing through 5 convolutional layers, the original input data, which has the shape of $256 \times 256 \times 4$, is reduced to the shape of $6 \times 6 \times 256$ and then flattened to enter the linear fully connected layers of $2048$, $2048$, and $1024$ neurons, respectively. In these layers, as well as in the convolution layer, the ReLU nonlinear activation function is utilized. Finally, the output layer consists of 3 neurons with the identity activation function. The architecture of the surrogate is described in Table \ref{NN_architecture}.

\begin{table}
\centering
\caption{The Architecture of the surrogate}
\label{NN_architecture}
\begin{tabular}[t]{lcc}
\hline
Layer & Type & Output size 
\\
\hline
 &  Input  &$256 \times 256 \times 4$ \\
1 & Conv2D & $127 \times 127 \times 24$ \\
2 & Conv2D & $62 \times 62 \times 48$ \\
3 & Conv2D & $30 \times 30 \times 96$ \\
4 & Conv2D & $14 \times 14 \times 192$ \\
5 & Conv2D & $6 \times 6 \times 256$ \\
7 & Linear & $2048$ \\
8 & Linear & $2048$ \\
9 & Linear & $1024$\\
10 & LinearOutput & $3$  \\
\hline
\end{tabular}
\end{table}%
\FloatBarrier

\subsection{Metrics for surrogate accuracy assessment}
Let ${\mathcal{D} = \{(\vc{X}_j, \vc{K}_j)\}_{j=1}^{D}}$ be a dataset of i.i.d. samples, where $\vc X_j \in \mathbb{R}^{4 \times F}$ is a vector of input features (with $F$ representing the number of pixels) and $\vc{K}_j \in \mathbb{R}^3$ is the vector of components $k_{xx}$, $k_{xy}$, and $k_{yy}$ of the corresponding equivalent hydraulic conductivity tensor.
Since supervised learning is carried out, $\mathcal{D}$ is divided into learning (=training) samples $\mathcal{L}$, validation samples $\mathcal{V}$, and test samples $\mathcal{T}$.  
Given $\mathcal{L}$, the~surrogate learns a function $f(\vc{x}): \mathbb{R}^{|\vc X|} \rightarrow \mathbb{R}^3$. 
During the training process, the accuracy of the surrogate's predictions is evaluated by a loss function, in our case, the mean squared error (MSE): $\gamma(\mathcal M, f) =  \frac{1}{|\mathcal M|}\|\vc{K}_{\mathcal{M}} - f(\vc{X}_{\mathcal{M}})\|^2_2$, where $\mathcal{M} \subseteq \mathcal D$.
The training loss $\gamma(\mathcal{L}, f)$ and validation loss $\gamma(\mathcal V, f)$ are observed. 

For the sake of comparing the accuracy of surrogates, we utilize:
 \begin{itemize}
     \item the coefficient of determination for component $i$ of $\vc{K}$:
\begin{equation}
R^2_i = 1 - \frac{\sum_{\mathcal M}{(\vc{K}^i - f(\vc{X}))^2}}{\sum_{\mathcal M}{(\vc{K}^i - \overline{\vc{K}^i_{\mathcal{M}}})^2}},   
\end{equation}
where $\overline{\vc{K}^i_{\mathcal{M}}}$ is the average of $\vc K^{i}$ over dataset $\mathcal M$. For further analysis, we also employ $\overline{R^2} = \frac{R^2_{xx} + R^2_{xy} + R^2_{yy}}{3}$. The closer $R^2$ is to $1$, the better.
\item normalized root mean squared error for component $i$ of $\vc{K}$: 
\begin{equation}
{N\!R\!M\!S\!E}_{i} = \frac{\sqrt{\frac{1}{|\mathcal M|}\|\vc{K}^i_{\mathcal{M}} - f(\vc{X}_{\mathcal{M}})\|^2_2}}{\std(\vc{K}^i_{\mathcal{M}})},
\end{equation}
where $\std(\vc{K}^i_{\mathcal{M}})$ is the standard deviation of $\vc{K}^i$ over the dataset $\mathcal{M}$.
The lower the ${N\!R\!M\!S\!E}$, the better. Moreover, if the ${N\!R\!M\!S\!E}$ is above $1$, we can use a simple random generator instead of a complex~surrogate.
\end{itemize}

\section{Dataset configuration}\label{dataset_config}
Datasets are constructed following our homogenization procedure described in Section \ref{num_hom_multiscale}.
The original domain $\Omega_{o} = (0, 100)^2$ is extended to the domain $\Omega = (0, 114.28)^2$. A DFN is generated with the density $\rho^{\prime}_{2D} = 10.0$, comprising approximately $1500$ fractures. Fracture centers are placed in $\Omega$ by the Poisson process. We consider a power-law distribution of fracture lengths ranging from $\rmin = 4.325$ to $\rmax = 100.0$, and the orientation angle is uniformly distributed. Based on this, the geometry with a mesh with step $h=4.325$ is generated.
The extended domain $\Omega$ is covered by $225$ overlapping homogenization blocks of size $(14.28, 14.28)$, where the centers of corner homogenization blocks are placed in the corners of the original domain $\Omega_{o}$. 
The matrix SRF is generated with parameters $\vmu = [-6.  -5.8]$,
$
  \vc \Sigma=
  \begin{bmatrix}
    0.25 & 0.2\\
    0.2  & 0.25
  \end{bmatrix}
$.

We have analyzed the problem for diverse fracture-to-matrix hydraulic conductivity ratios $K_f/K_m$, where \\ ${K_f/K_m \in \{\num{1e3}, \num{1e5}, \num{1e7}\}}$ was selected as representative enough for this study. 
Consequently, three datasets with different $K_f/K_m$ ratios of inputs were formed:
\begin{itemize}
    \item \textbf{Dataset $\mathcal{A}$}: $K_f/ K_m = \num{1e3}$,
    \item \textbf{Dataset $\mathcal{B}$}: $K_f/ K_m = \num{1e5}$,
    \item \textbf{Dataset $\mathcal{C}$}: $K_f/ K_m = \num{1e7}$.
\end{itemize}
Each dataset contains $75,000$ samples. In order to provide predictions for a range of matrix SRF correlation lengths up to $\lambda = 25$, correlation lengths $\lambda \in \{0, 10, 25\}$ are equally represented in the datasets. 
The $75,000$ samples are split into $80\%$ training samples and $20\%$ test samples.
Within the training samples, $80\%$ are allocated for actual training, while $20\%$ are reserved for validation.

\begin{figure*}
  \centering
  \begin{subfigure}{\textwidth}
    \includegraphics[width=0.32\linewidth]{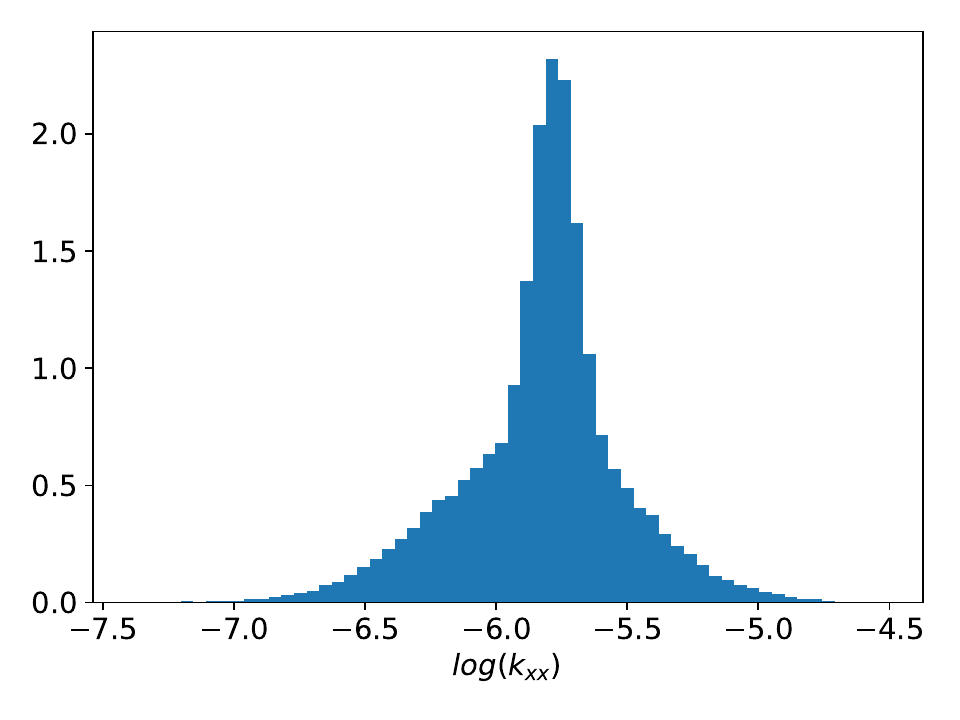}
    \includegraphics[width=0.32\linewidth]{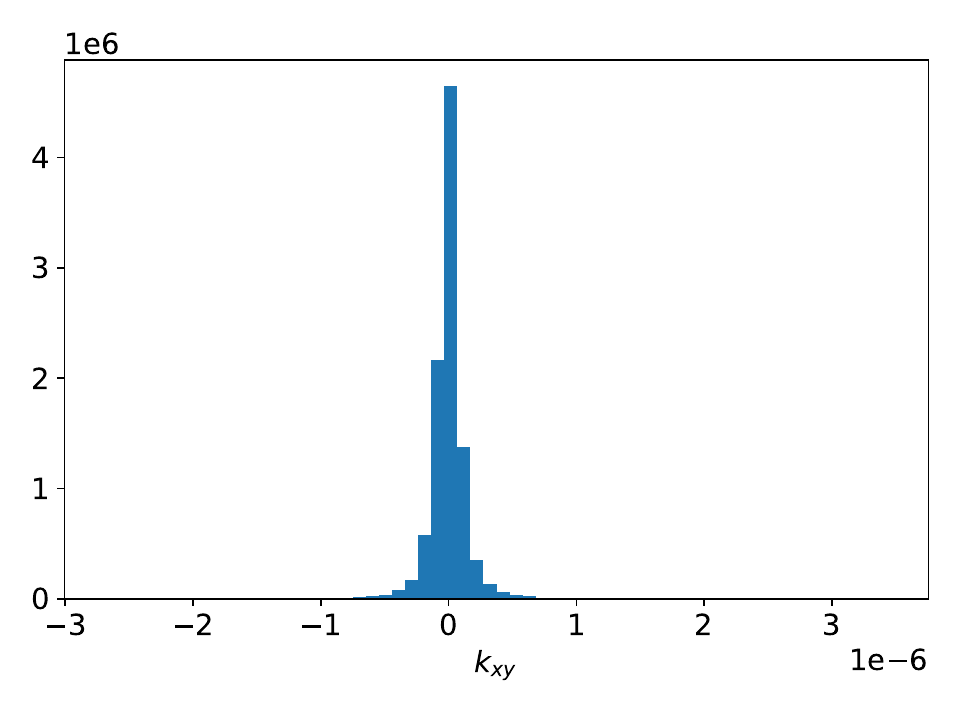}
    \includegraphics[width=0.32\linewidth]{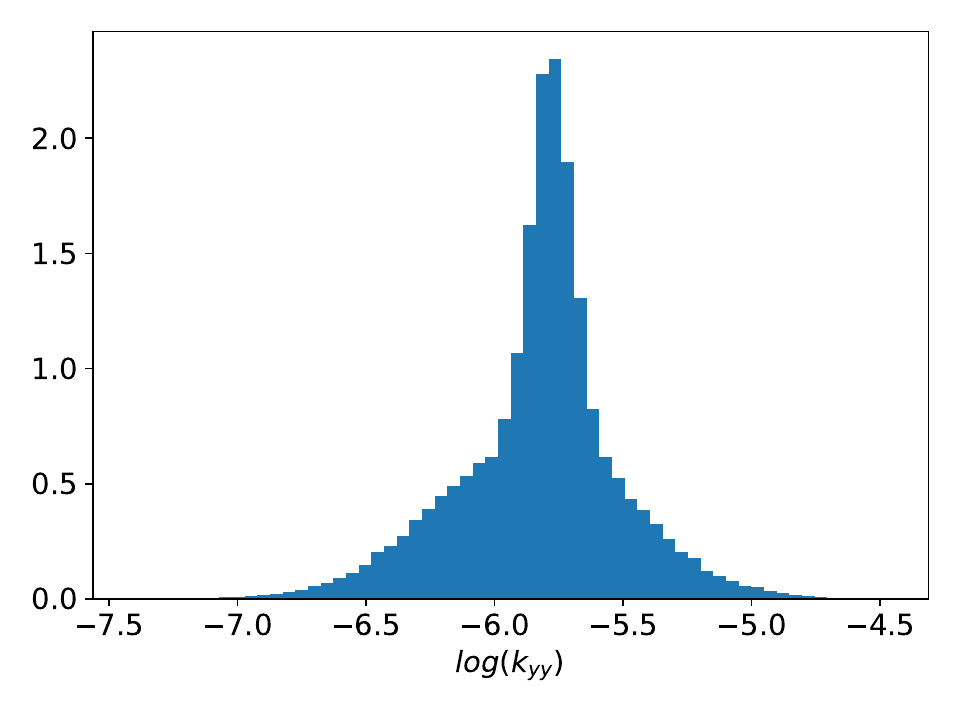}
    \caption*{Dataset $\mathcal{A}$}
  \end{subfigure}
  \centering
  \medskip

  \begin{subfigure}{\textwidth}
    \includegraphics[width=0.32\linewidth]{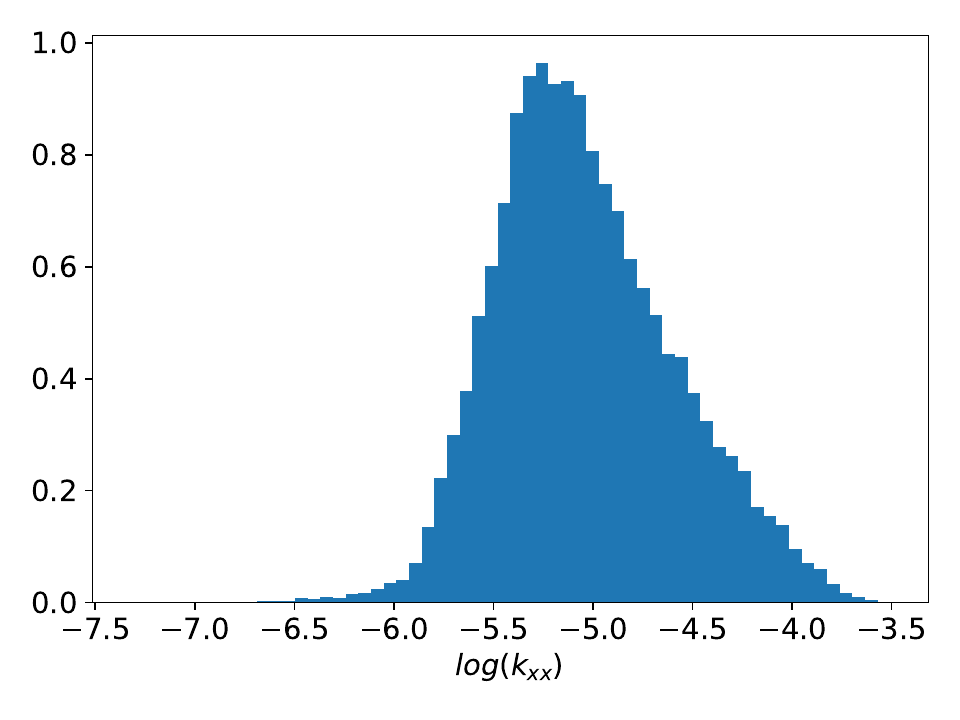}
    \includegraphics[width=0.32\linewidth]{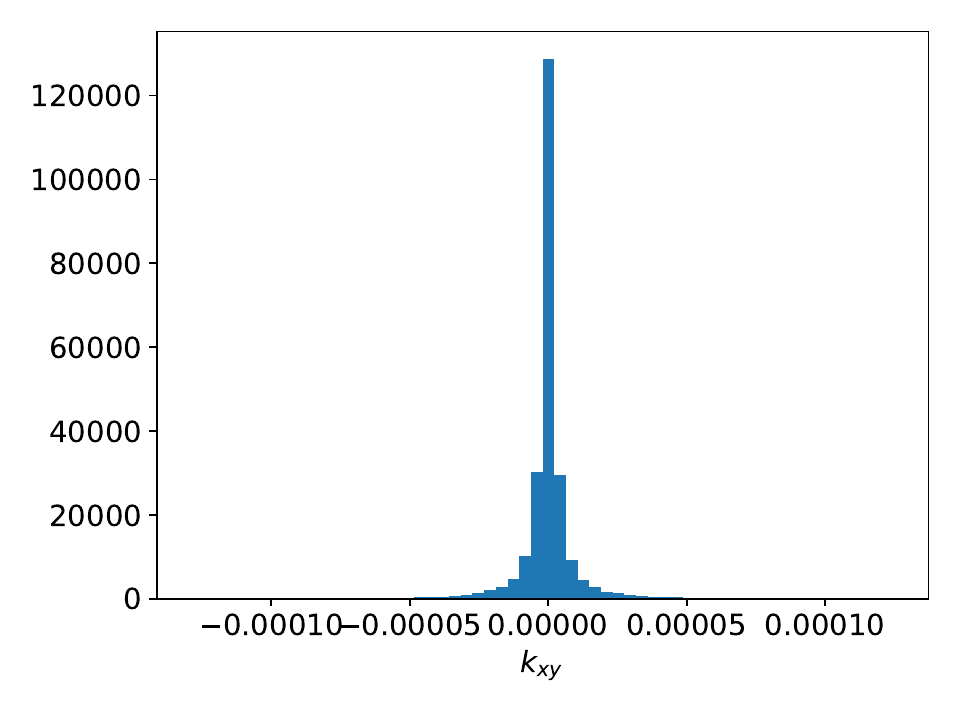}
    \includegraphics[width=0.32\linewidth]{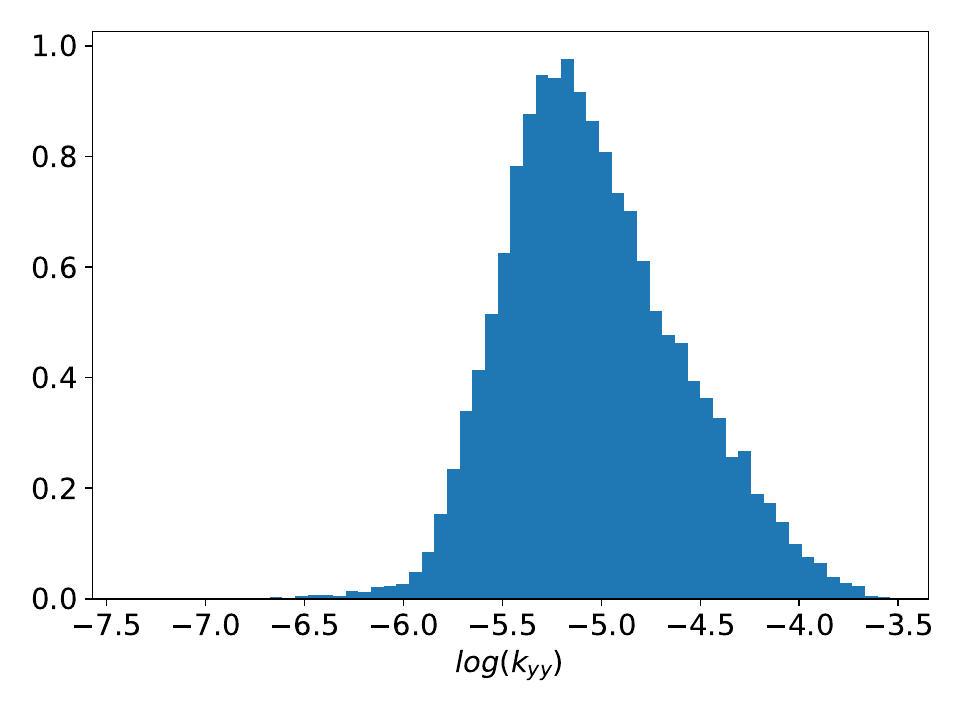}
      \caption*{Dataset $\mathcal{B}$}
  \end{subfigure}

  \medskip

  \begin{subfigure}{\textwidth}
    \includegraphics[width=0.32\linewidth]{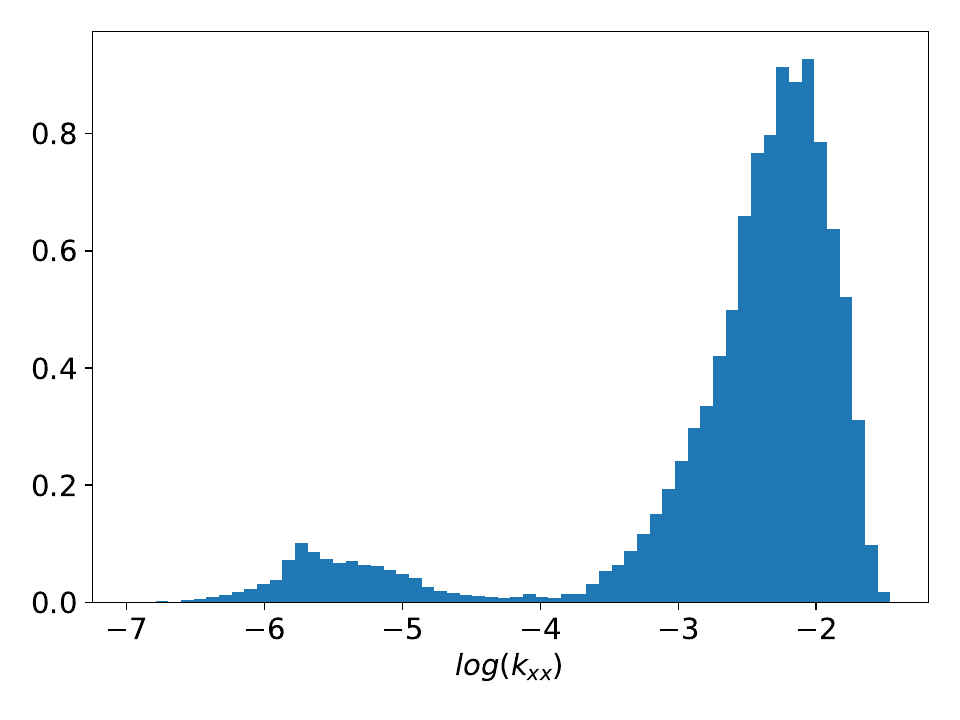}
    \includegraphics[width=0.32\linewidth]{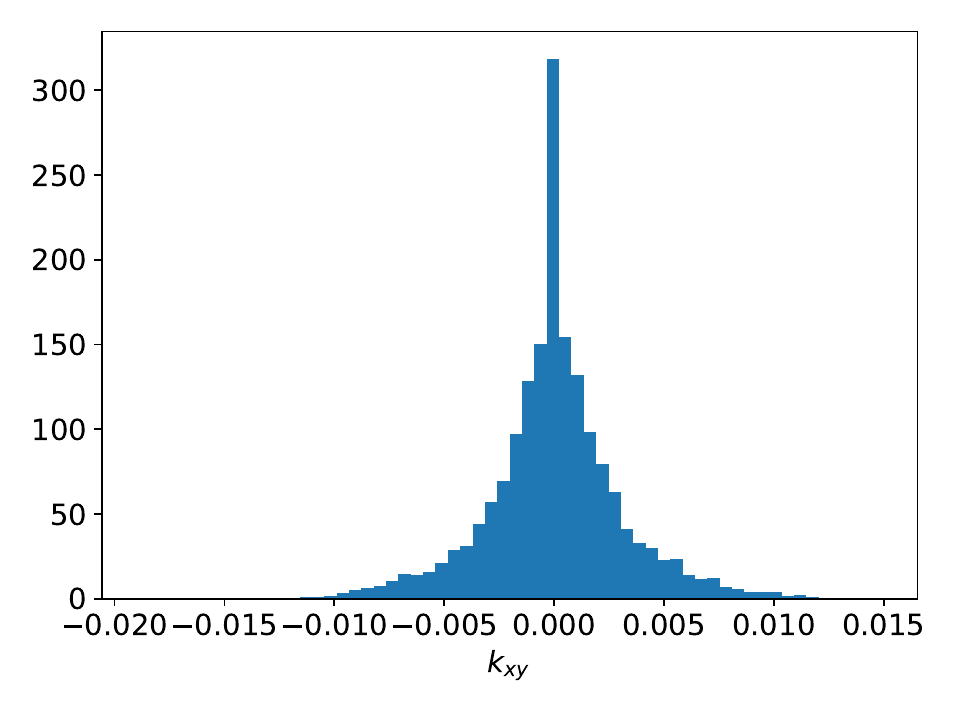}
    \includegraphics[width=0.32\linewidth]{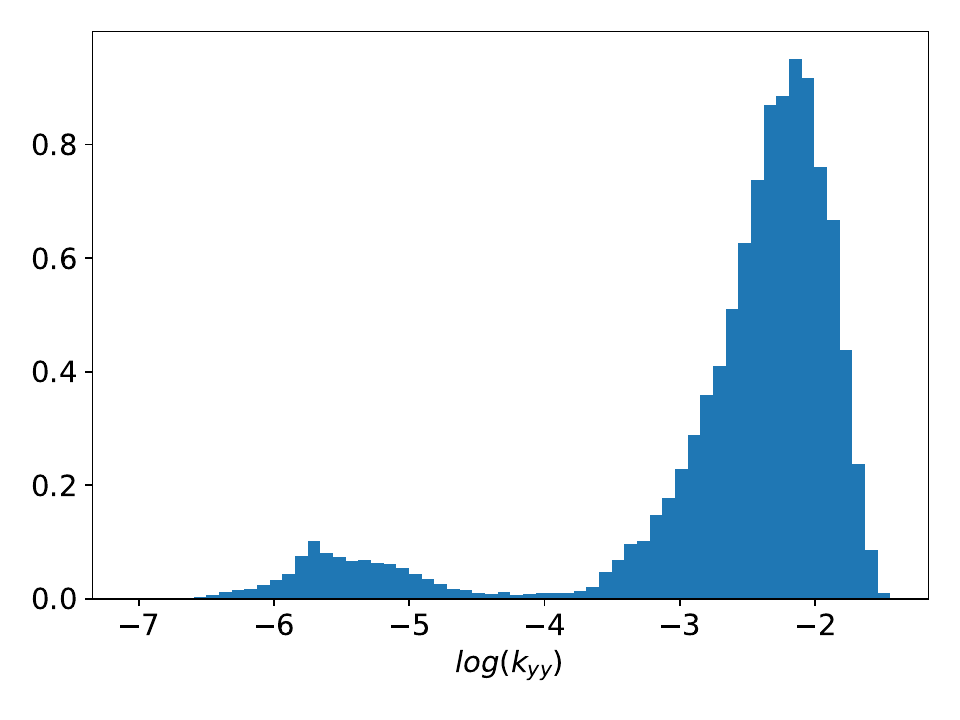}
      \caption*{Dataset $\mathcal{C}$}
  \end{subfigure}
  \caption[Distributions of $\tn K^{eq}$ components for datasets of different $K_f/K_m$]{Distributions of $\tn K^{eq}$ components for datasets of different $K_f/K_m$, Dataset $\mathcal{A}$ for $K_f/ K_m = \num{1e3}$, Dataset $\mathcal{B}$ for $K_f/ K_m = \num{1e5}$, and Dataset $\mathcal{C}$ for $K_f/ K_m = \num{1e7}$.}
  \label{fig:datasets}
\end{figure*}

Distributions of equivalent hydraulic conductivity tensors for datasets $\mathcal{A}$,$\mathcal{B}$, and $\mathcal{C}$ are depicted in Fig.\ref{fig:datasets}. Considering $k_{xx}$ and $k_{yy}$, we observe that for Dataset $\mathcal{A}$, the equivalent tensor is predominantly determined by matrix elements. For Dataset $\mathcal{B}$, we note an increased impact of fractures. However, there is still no visible percolation threshold characterized by a sharp increase in equivalent hydraulic conductivity at the moment of significant interconnection of fractures. This phenomenon has been observed in studies such as \cite{HE2023592, Lang20140815}. 
As the ratio $K_f/ K_m$ increases, the gap between non-percolated and percolated samples widens. We observe orders of magnitude difference for Dataset $\mathcal{C}$. Regarding $k_{xy}$, the shape of the data is very similar across datasets.

The equivalent tensor distributions exhibit significant dissimilarities in scale and shape. Consequently, training a surrogate model on the combined dataset of all presented $K_f/ K_m$ ratios has proven to be extremely challenging and has resulted in unsatisfactory accuracy thus far. Therefore, we resort to training separate surrogate models: Surrogate A for Dataset $\mathcal{A}$, Surrogate B for Dataset $\mathcal{B}$, and Surrogate C for Dataset $\mathcal{C}$. Subsequently, based on the $K_f/ K_m$ ratio of a particular input $\vc X_j$, the most suitable surrogate is selected for prediction.

\subsection{Dataset preprocessing}
In order to facilitate the surrogate's training, hydraulic conductivity tensors in the dataset inputs $\vc X$ and outputs $\vc K$ are preprocessed in the same way.
For each sample $j$, the average matrix hydraulic conductivity $\overline{\vc X^j}_{m}$ (average over the entire tensor) is calculated.

First, inputs and outputs are normalized so that hydraulic conductivities of different scales with the same $K_f/K_m$ value can be accommodated:
\begin{equation}
\begin{aligned} 
\vc X^j &= \vc X^j / \overline{\vc X^j}_{m}\\
\vc K_j &= \vc K_j / \overline{\vc X^j}_{m}
\end{aligned}
\end{equation}
Second, data is independently standardized for input and output hydraulic conductivities, where mean ($\text{avg}$) and standard deviation ($\text{std}$) are calculated based on the training set $\mathcal{L}$. Each component of a hydraulic conductivity tensor is preprocessed as follows:
$$
\tilde{k}^j_{xx} = \frac{\log(k^j_{xx}) - \text{avg}_{\mathcal{L}}(\log(k_{xx}))}{\text{std}_{\mathcal{L}}(\log(k_{xx}))},
$$
$$
\tilde{k}^j_{xy} = \frac{k^j_{xy} - \text{avg}_{\mathcal{L}}(k_{xy})}{\text{std}_{\mathcal{L}}(k_{xy})},
$$
$$
\tilde{k}^j_{yy} = \frac{\log(k^j_{yy}) - \text{avg}_{\mathcal{L}}(\log(k_{yy}))}{\text{std}_{\mathcal{L}}(\log(k_{yy}))}.
$$

\section{Results}\label{results_section}
In this section, we present the prediction accuracy of trained surrogates and discuss their prediction capabilities on datasets with different matrix SRF correlation lengths $\lambda$ and fracture densities $\rho^{\prime}_{2D}$.
Additionally, we employ trained surrogates for upscaling a DFM model within two macroscale problems to demonstrate the impact of prediction accuracy of $\tn K^{eq}$ on a coarser scale.
We use $125$ epochs for surrogate training, and the model with the lowest validation loss is saved for further use. During training, the initial learning rate $\alpha = 0.0025$ adaptively decreases to 10 percent of the previous value when the validation loss does not decrease in 10 subsequent epochs.

\begin{figure*}
  \centering
  \begin{subfigure}{\textwidth}
    \includegraphics[width=0.32\linewidth]{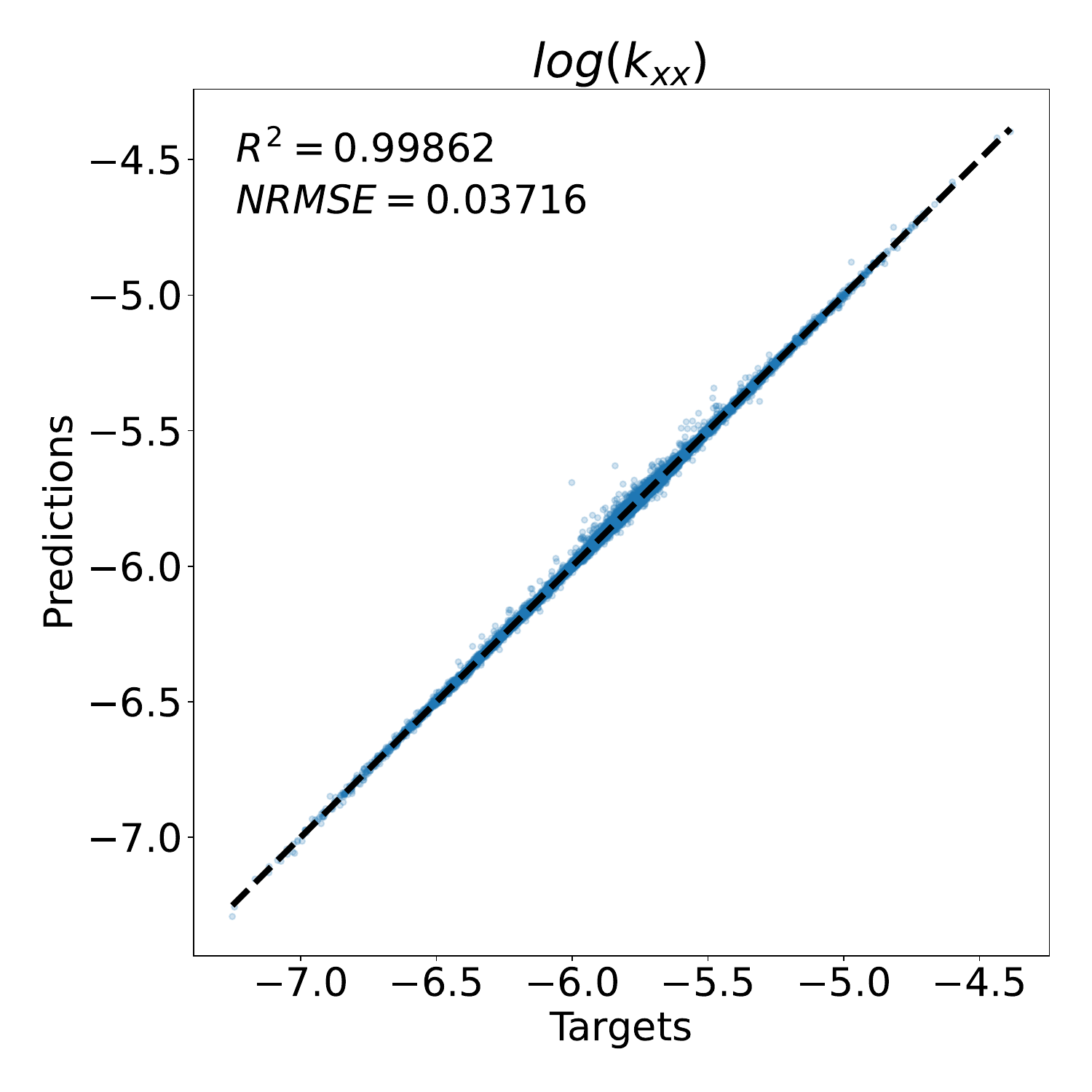}
    \includegraphics[width=0.32\linewidth]{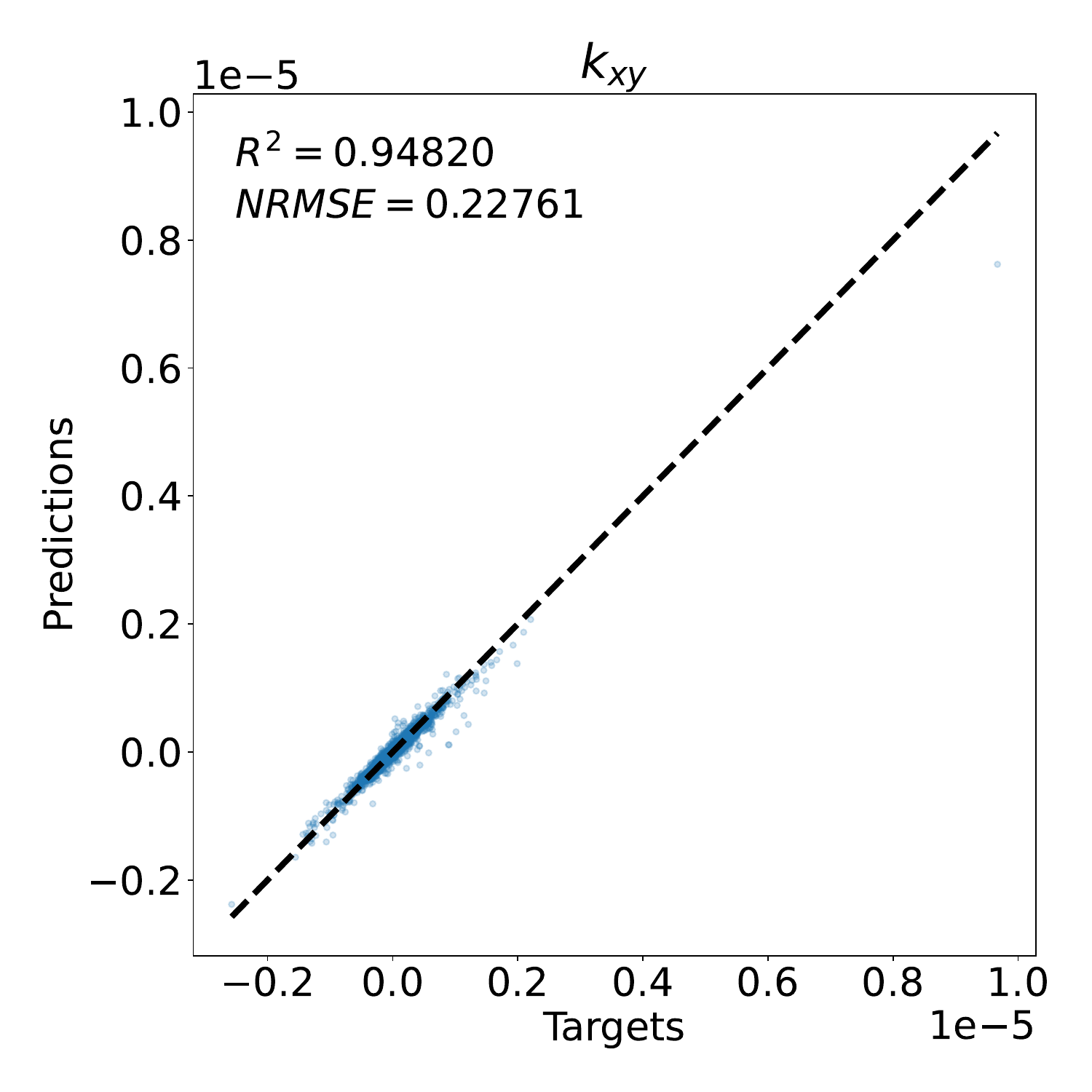}
    \includegraphics[width=0.32\linewidth]{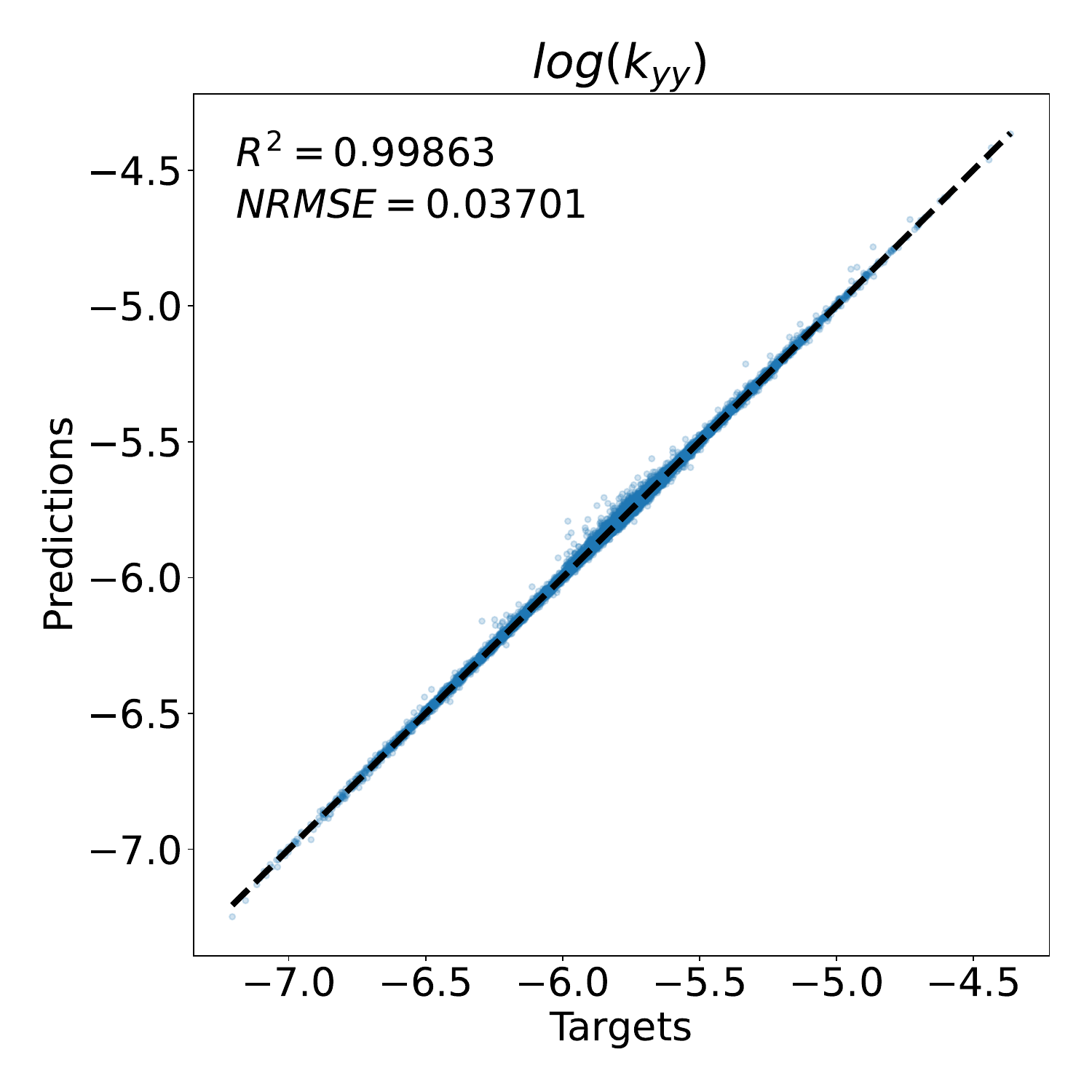}
    \captionsetup{justification=centering}
    \caption*{Surrogate A}
  \end{subfigure}
  \centering
  \medskip
  \begin{subfigure}{\textwidth}
    \includegraphics[width=0.32\linewidth]{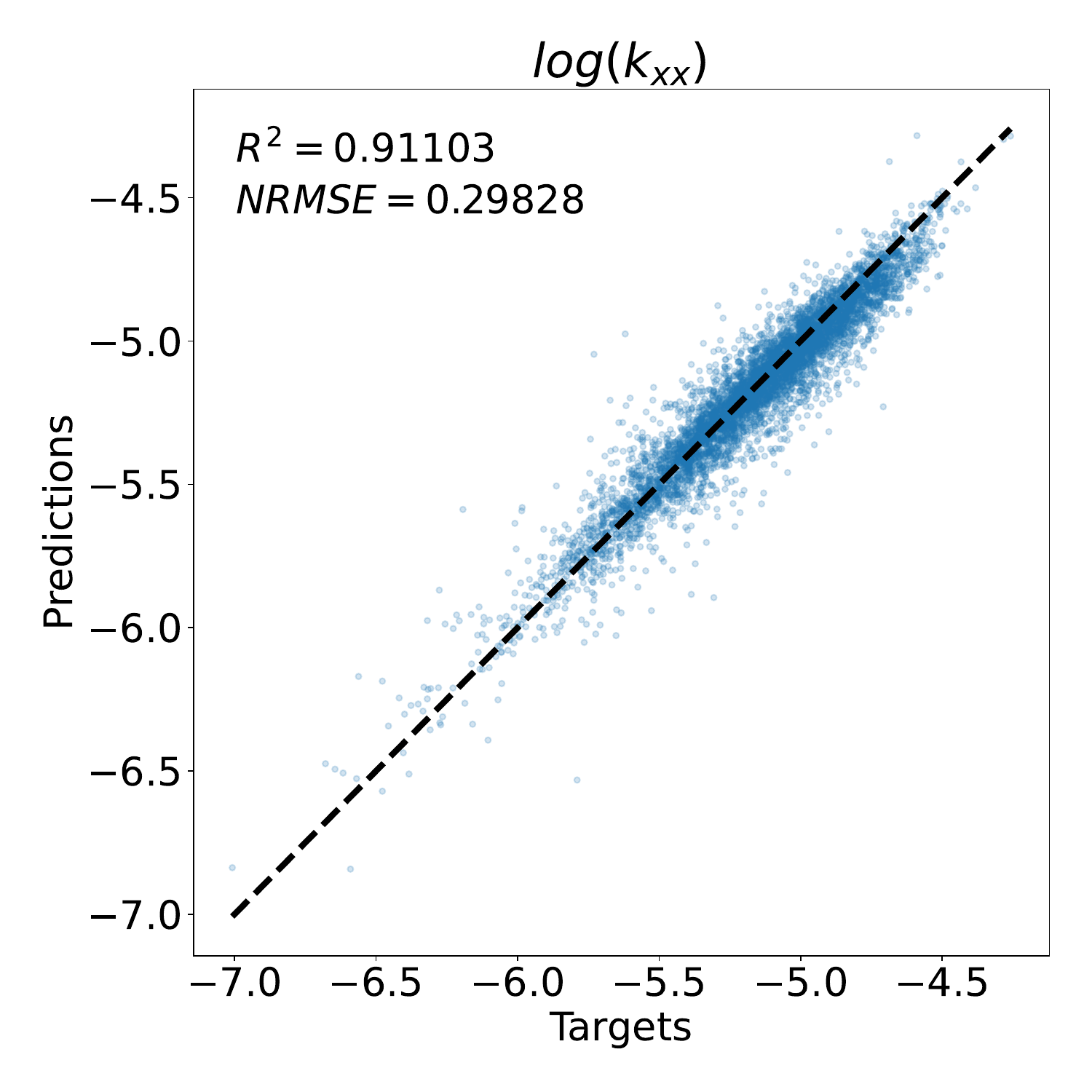}
    \includegraphics[width=0.32\linewidth]{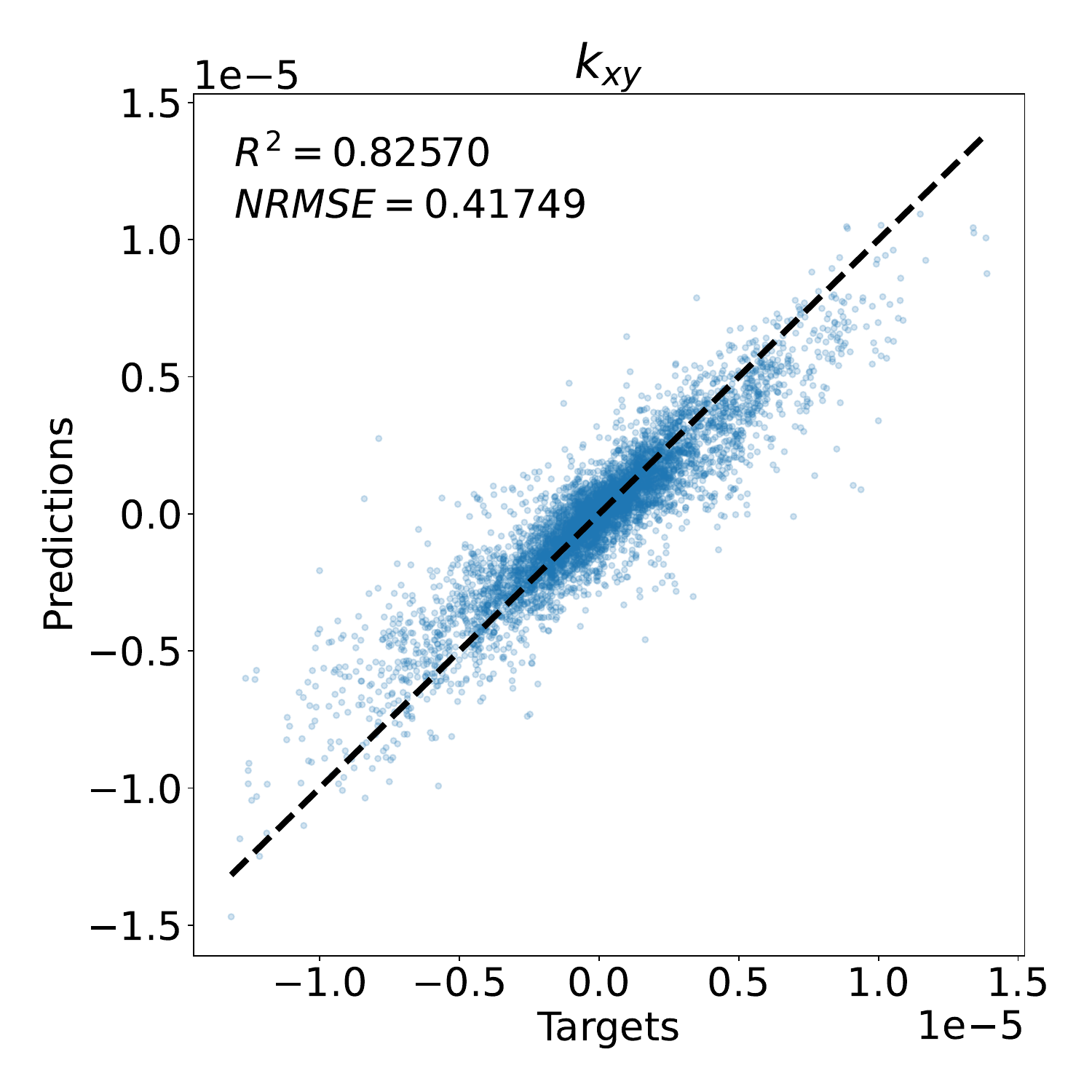}
    \includegraphics[width=0.32\linewidth]{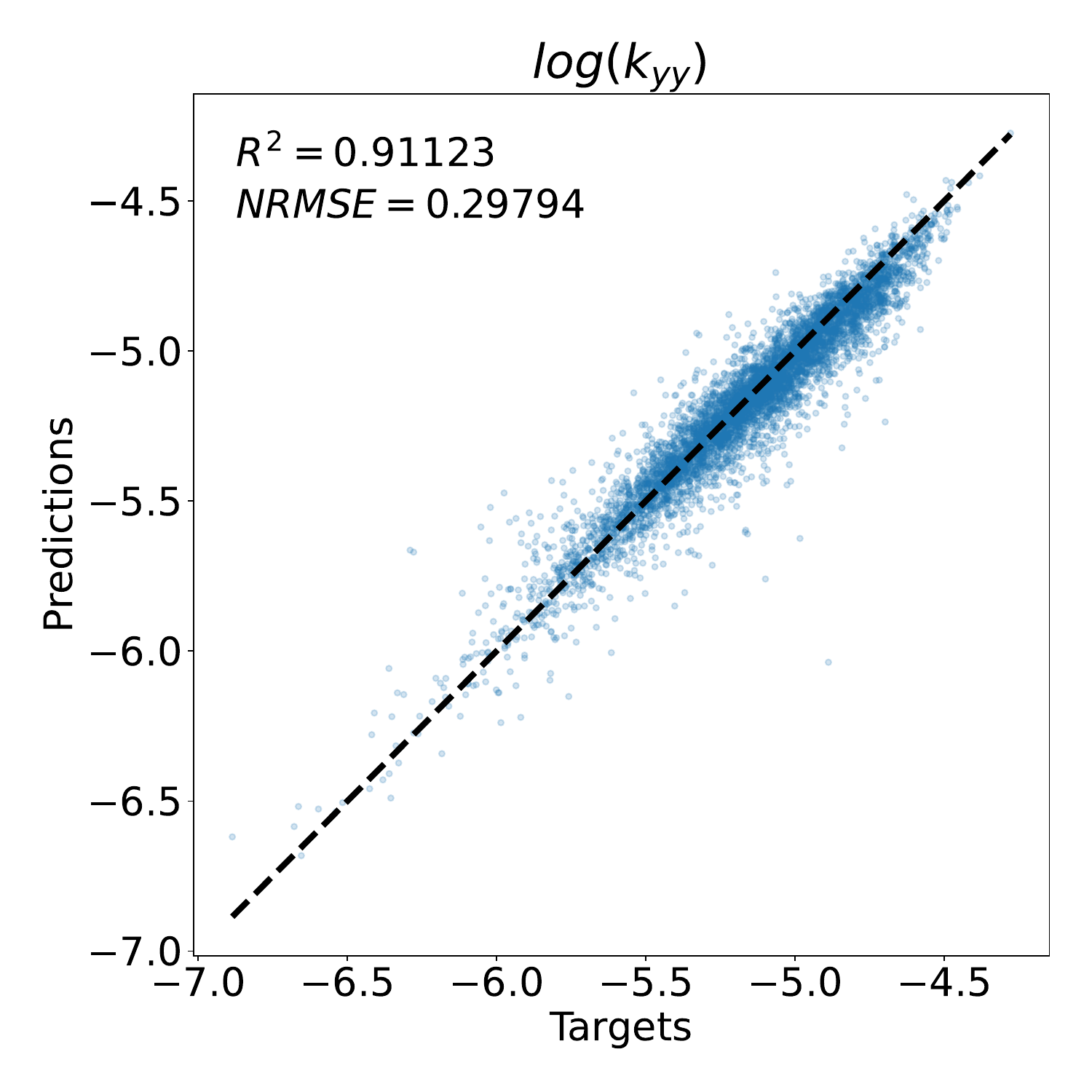}
      \caption*{Surrogate B}
  \end{subfigure}
  \medskip
  \begin{subfigure}{\textwidth}
    \includegraphics[width=0.32\linewidth]{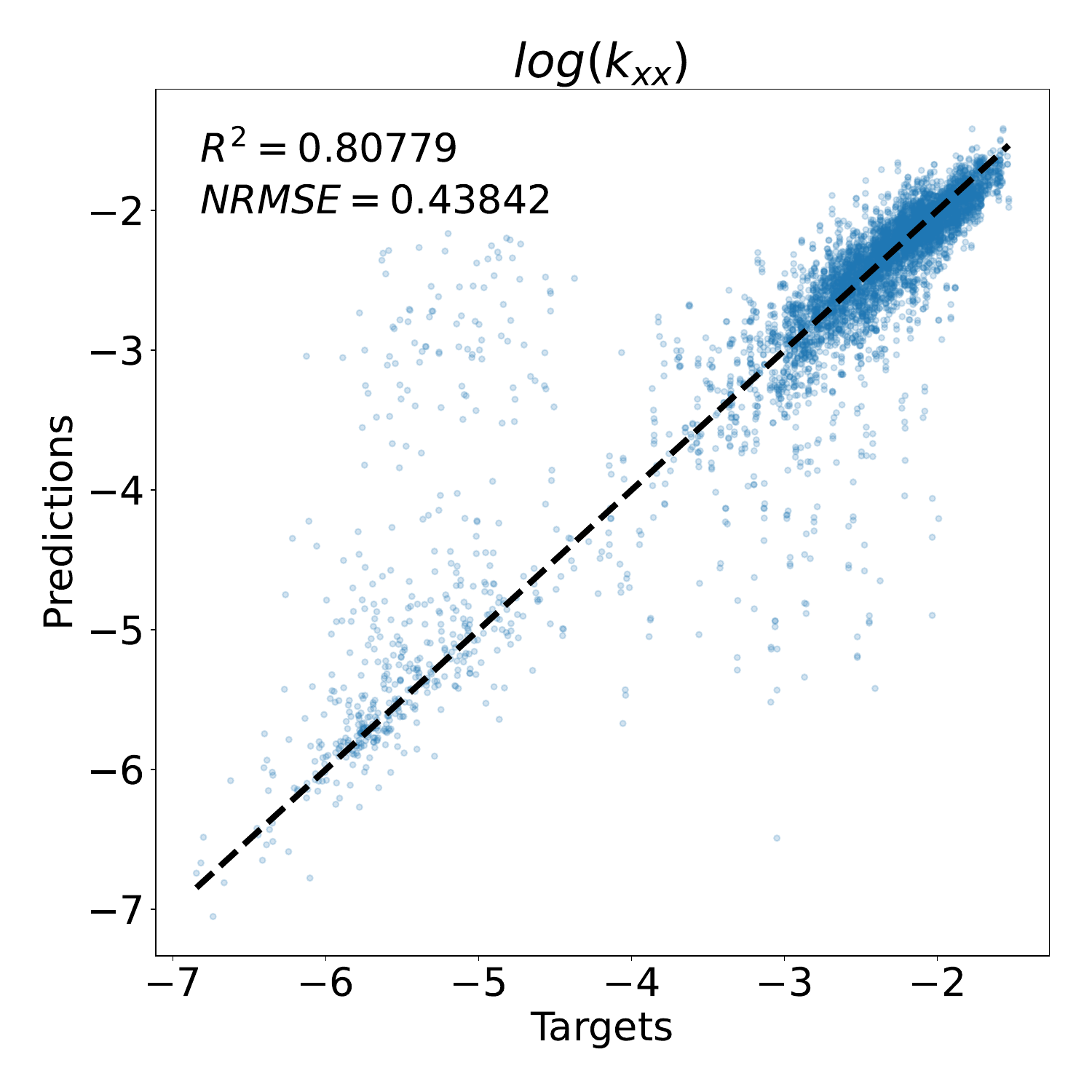}
    \includegraphics[width=0.32\linewidth]{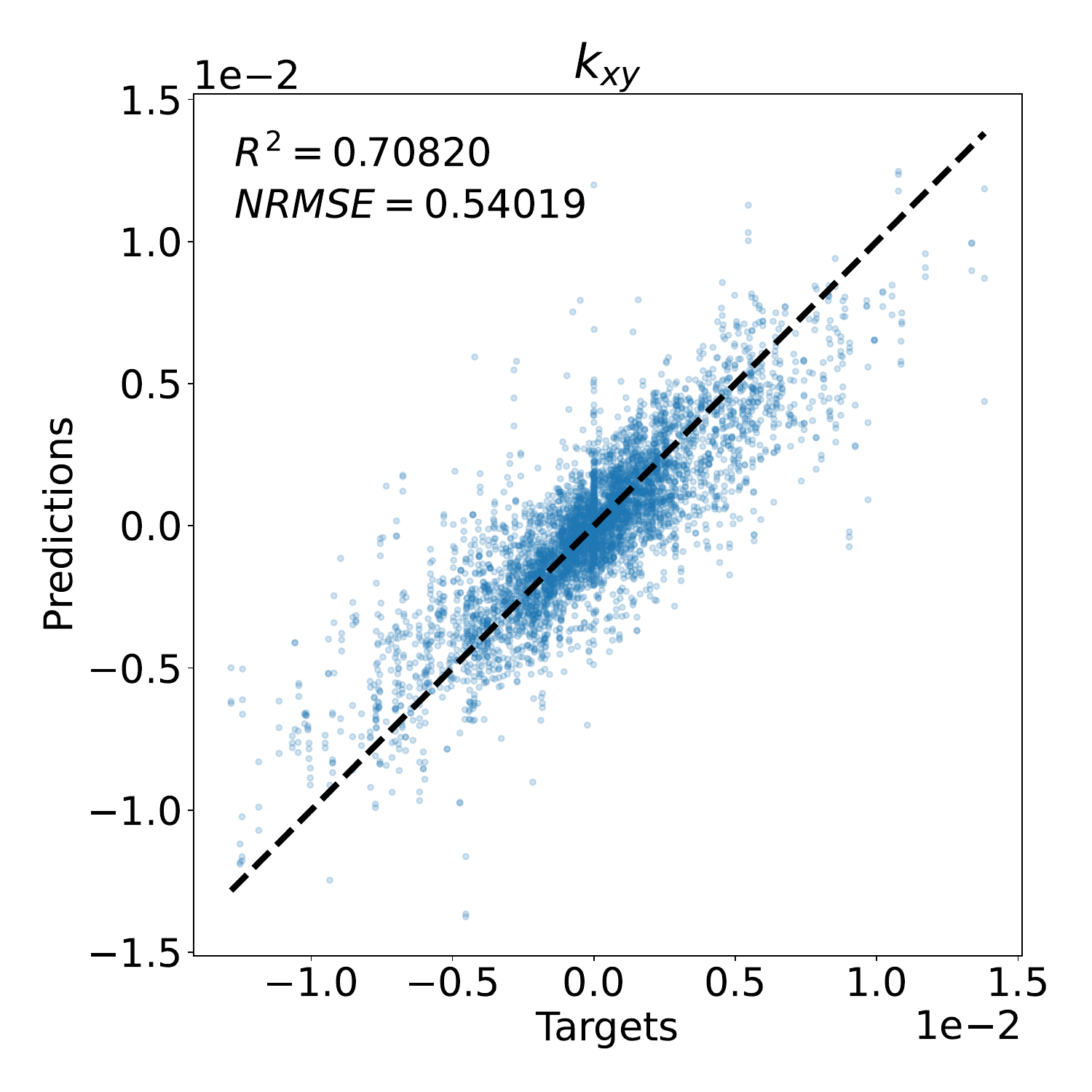}
    \includegraphics[width=0.32\linewidth]{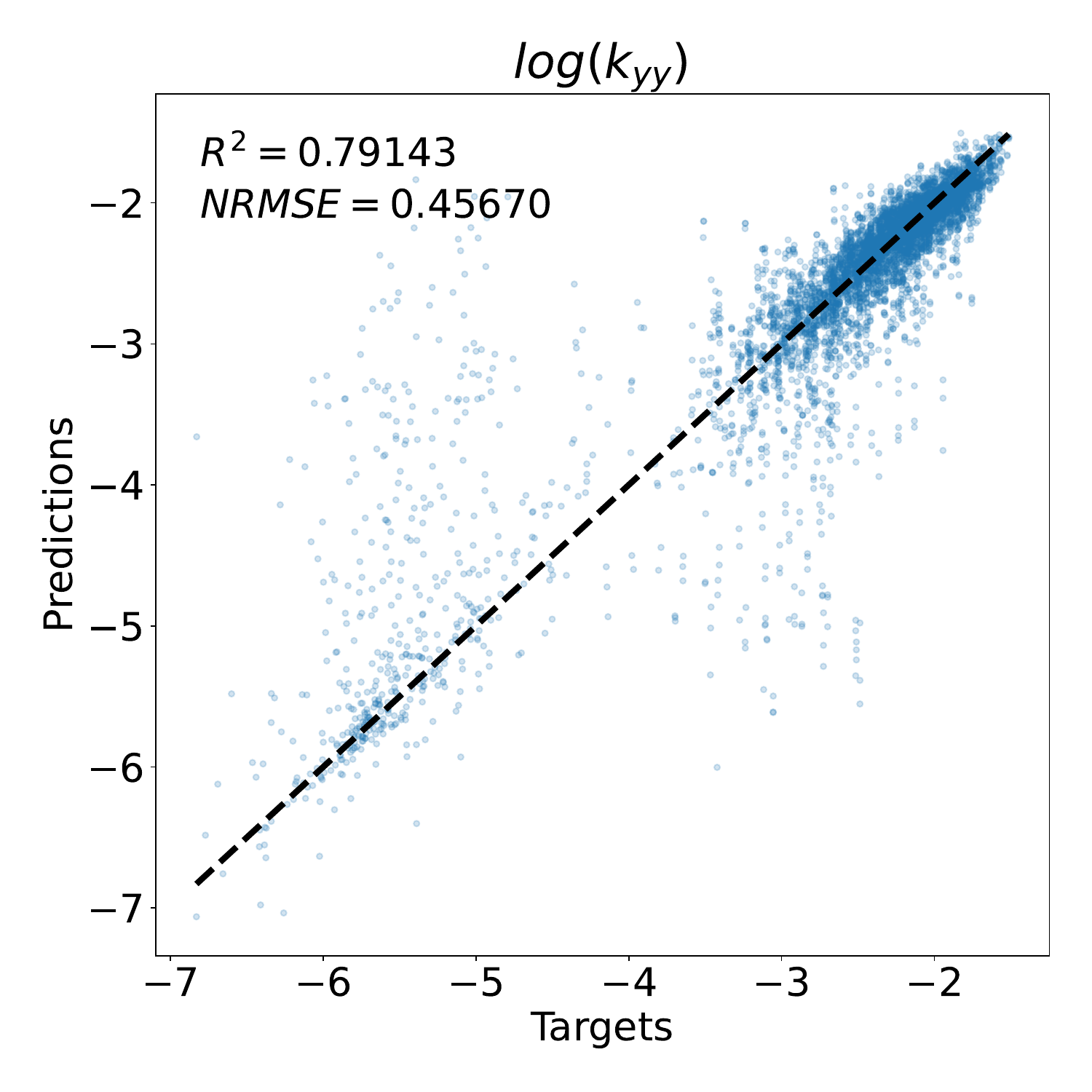}
      \caption*{Surrogate C}
  \end{subfigure}
  \caption[Accuracy of the surrogates on test datasets]{Accuracy of the surrogates on test datasets}
  \label{fig:surrogates_accuracy}
\end{figure*}

Fig.\ref{fig:surrogates_accuracy} shows the targets-predictions plot along with $R^2$ and NRMSE values for the test sets. The best results ($R^2 > 0.95$) were achieved for Surrogate A, where the distribution of tensor components is relatively narrow. 
These results are comparable to or better than those provided by Pal et al. \cite{pr11020601}, who do not consider fractures, or by Zhu et al. \cite{ZHU2023111186}, who only consider DFN models, meaning the impact of the matrix on the equivalent hydraulic conductivity is not considered.

As the training data distribution becomes broader and more complex (Fig. \ref{fig:datasets}), the surrogate accuracy deteriorates (Fig. \ref{fig:surrogates_accuracy}). 
This is particularly visible for Surrogate C, where the bimodal distributions of $\log(k_{xx})$ and $\log(k_{yy})$ lead to poor accuracy for samples between two distant peaks. This is due to the lack of this specific data in the training set compared to the data in the distribution peaks (Fig. \ref{fig:datasets}), especially the larger peak, where the predictions tend to be more accurate. This effect intensifies with increasing $K_f/K_m$.
For all surrogates, the off-diagonal component of $\tn K^{eq}$ is predicted with higher error than the diagonal components. Due to the nature of $k_{xy}$, we cannot carry out logarithmic-like preprocessing. Consequently, $k_{xy}$ suffers from a scarcity of training data in tail areas, where a more pronounced dispersion of predictions is observed.

Although using a larger training dataset generally improves the prediction accuracy of deep neural networks, it may not effectively address the issues associated with underrepresented data if the distribution of the additional training data remains the same.
Instead of simply increasing the size of the training set, we argue that to enhance prediction accuracy, one should focus on generating more data in the tails for Surrogate B or in the area between the lower and higher peaks for Surrogate C. This process is challenging because we lack a straightforward method to generate DFM homogenization blocks for a given equivalent hydraulic conductivity tensor.

In the following subsections, we investigate the prediction accuracy of trained surrogates on datasets with different fracture densities $\rho^{\prime}_{2D}$ and on datasets of different correlation lengths $\lambda$ of matrix SRF. Subsequently, we discuss the speedup achieved by surrogates.

\subsubsection{Impact of fracture density}
\begin{figure}
 \centering
 \includegraphics[width=0.5\textwidth]{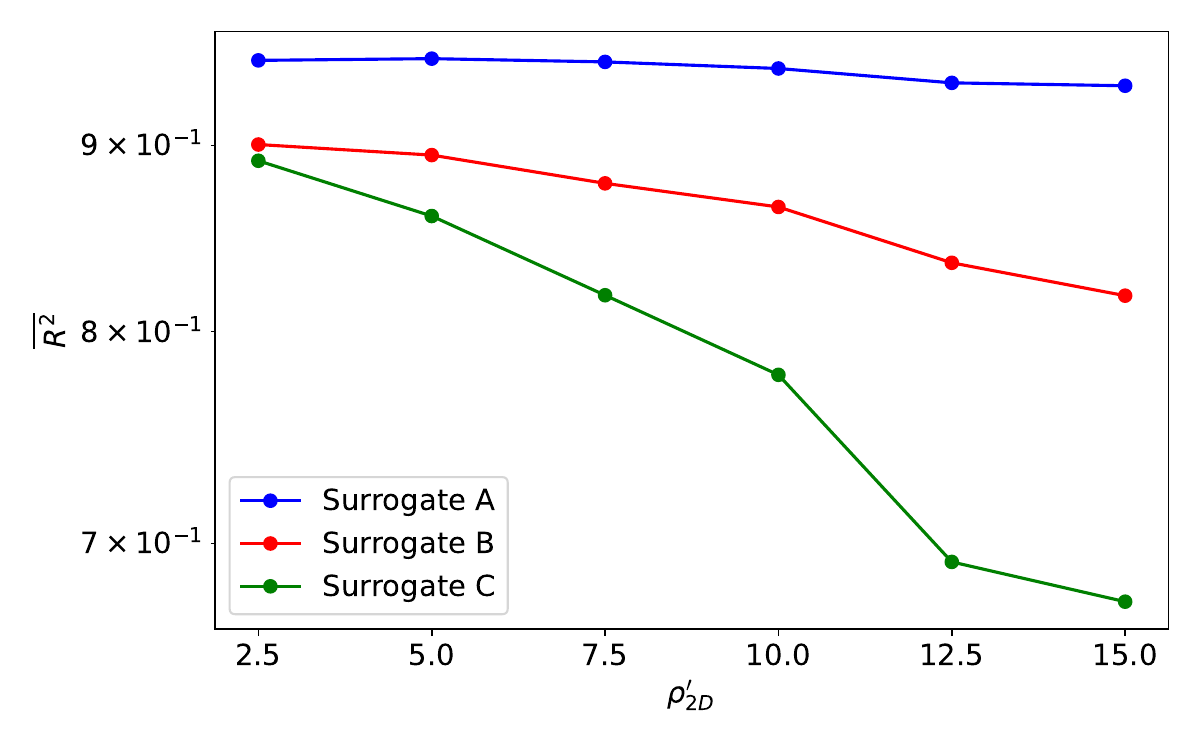}
  \caption{Impact of fracture density $\rho^{\prime}_{2D}$ on surrogate's accuracy}
 \label{fig:r2_rho}
\end{figure}

The impact of fracture density on surrogate accuracy has been previously discussed by Zhu et al. \cite{ZHU2023111186}, where they predicted equivalent permeability tensors based on DFN models. They generally observed a deterioration in prediction accuracy with increasing fracture density, which aligns with our findings.
Fig. \ref{fig:r2_rho} illustrates prediction accuracy for $\rho^{\prime}_{2D} \in \{2.5, 5, 7.5, 10, 12.5, 15\}$. Surrogates trained with $\rho^{\prime}_{2D}=10$ achieve increasingly better accuracy when predicting on datasets with $\rho^{\prime}_{2D} < 10$. 
Since $\rho^{\prime}_{2D}=10$ is prescribed for the entire domain $\Omega$, individual homogenization blocks may contain a number of fractures that corresponds to $\rho^{\prime}_{2D} < 10$. 
However, for $\rho^{\prime}_{2D} > 10$, $\overline{R^2}$ diminishes, with the most significant reduction observed in the case of Surrogate C, where fractures play a dominant role in $\tn K^{eq}$ calculation. Still, $\overline{R^2} \approx 0.68$ might suffice for some applications. If not, surrogate training on higher $\rho^{\prime}_{2D}$ might be a viable solution. 
It's important to note that even for Surrogate C trained on $\rho^{\prime}_{2D}=10$ data, some $k_{xy}$ components are inaccurately predicted for $\rho^{\prime}_{2D} \in \{2.5, 5\}$. The reason is that the $k_{xy}$ values span $7$ orders of magnitude ($\num{1e-9}$ to $\num{1e-2}$). Given the preprocessing of $k_{xy}$ and loss function, larger $k_{xy}$ values are predicted more accurately. 
When $\rho^{\prime}_{2D}$ is small enough, there are instances where the off-diagonal component of the hydraulic conductivity tensor is among the smallest in the dataset (around $\num{1e-8}$ and lower). This leads to substantial relative errors in predictions and, consequently, potentially non-positive definite predicted $\tn K^{eq}$ tensors.
We do not encounter this issue for $K_f/ K_m \leq \num{1e5}$, which is the maximum ratio investigated in studies on equivalent permeability/hydraulic conductivity calculation, such as \cite{Andrianov2022UpscalingOT, ZHU2023111186, Lang20140815, https://doi.org/10.1029/2001WR000756}.

\subsubsection{Impact of correlation length}
\begin{figure}
 \centering
 \includegraphics[width=0.5\textwidth]{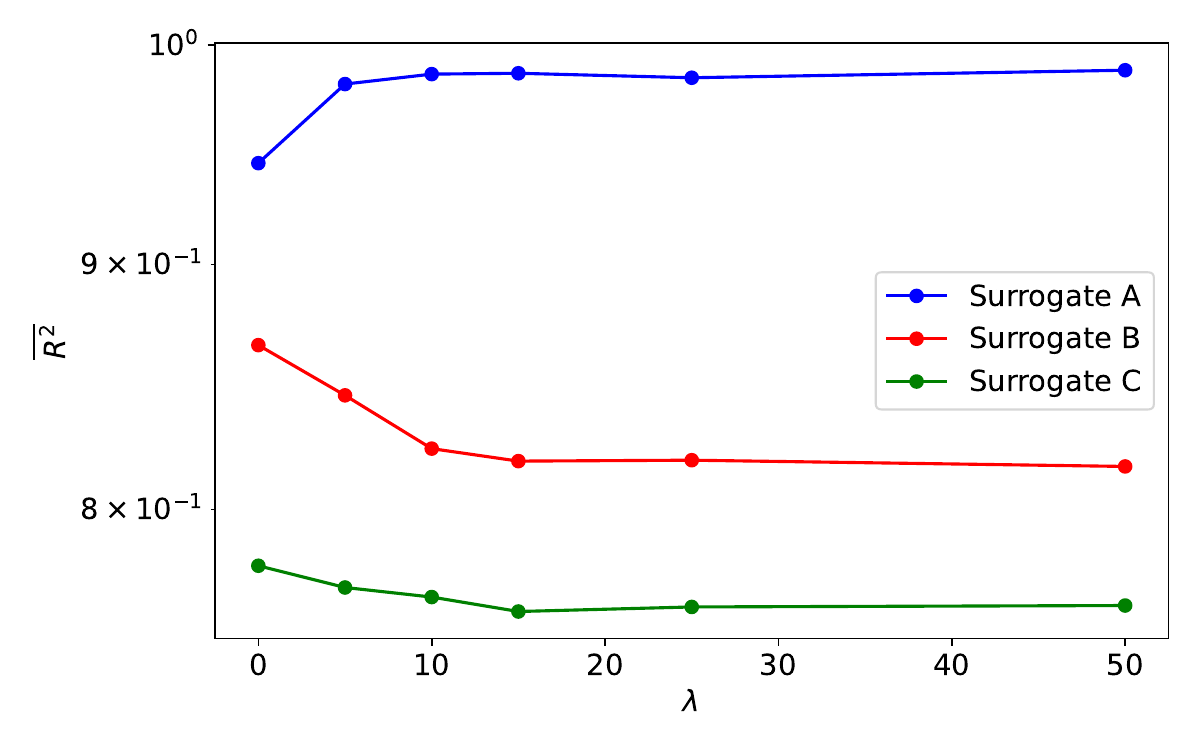}
 \caption{Impact of matrix correlation length $\lambda$ on surrogate's accuracy}
 \label{fig:r2_corr_len}
\end{figure}

We further investigate the impact of the matrix's SRF correlation length $\lambda$ on the prediction accuracy of our surrogates trained on $\lambda \in \{0, 10, 25\}$.
As shown in Fig. \ref{fig:r2_corr_len}, Surrogate A exhibits improved prediction accuracy with increasing $\lambda$. 
This trend is expected because a larger correlation length induces lower variability in the SRF field, making it easier for neural networks to learn. 
The opposite trend is observed in the cases of Surrogate B and Surrogate C. As mentioned earlier, these surrogates were trained on datasets with $K_f/K_m = \num{1e5}$ and $K_f/K_m = \num{1e7}$, where the hydraulic conductivity of fractures significantly influences $\tn K^{eq}$. Consequently, samples with very similar DFNs but dissimilar SRFs can result in similar $\tn K^{eq}$. As the correlation length increases, it becomes increasingly challenging for Surrogate B and Surrogate C to differentiate between these two cases, potentially leading to lower prediction accuracy. Despite the equal representation of all three correlation lengths $\lambda = \{0, 10, 25\}$ in the training dataset, Surrogate B and Surrogate C exhibit the best predictive performance for uncorrelated SRFs ($\lambda = 0$).

Another important observation is that although correlation lengths $\lambda = \{5, 15\}$ were not explicitly included in the training datasets, the corresponding predictive performance mirrors the trends observed for the included correlation lengths. 
This suggests that a wide range of correlation lengths on the homogenization block (of size $14.24 \times 14.28$) can be effectively covered by using datasets with $\lambda = \{0, 10, 25\}$, where the SRFs are generated following the procedure described in Section \ref{matrix_srf_section}.

\subsubsection{Computational cost reduction}
\begin{figure}
 \centering
 \includegraphics[width=0.5\textwidth]{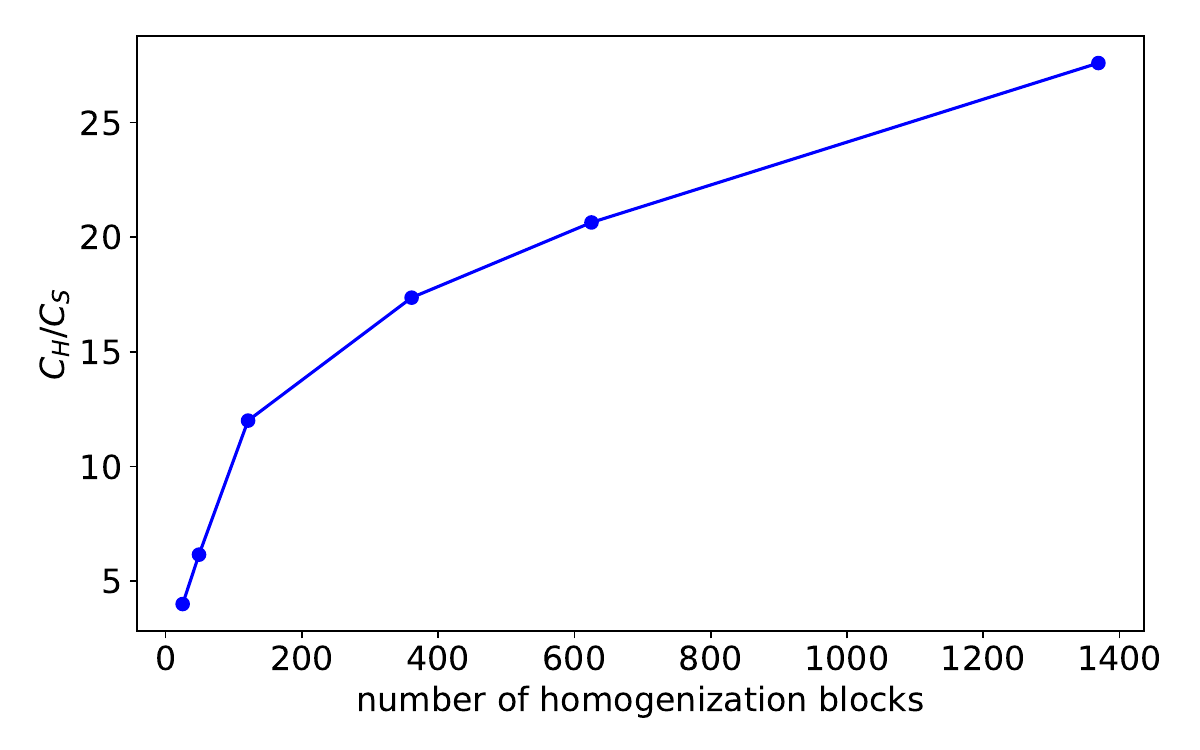}
 \caption{Speedup gained by surrogate predictions for different numbers of homogenization blocks}
 \label{fig:surrogate_speedup}
\end{figure}

We aim to utilize surrogates to speed up numerical homogenization. To this end, we compare the computational cost (CPU time) required for numerical homogenization ($C_{H}$) and surrogate prediction ($C_{S}$) across varying numbers of homogenization blocks covering the entire domain $\Omega$ (see Fig. \ref{fig:surrogate_speedup}). In the case of $C_{S}$, the predominant cost arises from the rasterization of an input mesh. The entire domain $\Omega$ is rasterized once and then divided into homogenization blocks. The prediction time of the neural network is negligible, and training time is excluded from consideration.
The cost of numerical homogenization $C_{H}$ encompasses the generation of a computational mesh and the execution of simulation software for each homogenization block.
Consequently, the acceleration achieved by using the surrogates increases with the number of homogenization blocks. For instance, with $25$ blocks, $C_{H}/C_{S} = 4$, while with $1369$ blocks, we attain $C_{H}/C_{S} \approx 28$. The latter scenario would correspond to the finest-scale samples, where we also assume the most accurate surrogate predictions.

\subsection{Upscaled equivalent hydraulic conductivity tensors in use}
We explore the impact of utilizing predicted equivalent tensors for the Aquifer problem and the Anisotropy problem, as described in Section \ref{sec:DFM}. In both cases, the effect of upscaling by numerical homogenization and by surrogates is compared.
We consider fine DFM samples with $\rho^{\prime}_{2D} = 10$, $\lambda=0$, $K_f/K_m \in \{\num{1e3}, \num{1e5}, \num{1e7}\}$. Trained Surrogate A, Surrogate B, and Surrogate C are used for predicting $\tn K^{eq}$. In accordance with the procedure explained in Section \ref{num_hom_multiscale}, homogenized/predicted $\tn K^{eq}$ are utilized for the matrix hydraulic conductivity of the upscaled DFM model.

\subsubsection{Aquifer problem}
\begin{figure}
  \centering
  \begin{subfigure}{0.4\textwidth}
    \includegraphics[width=\linewidth]{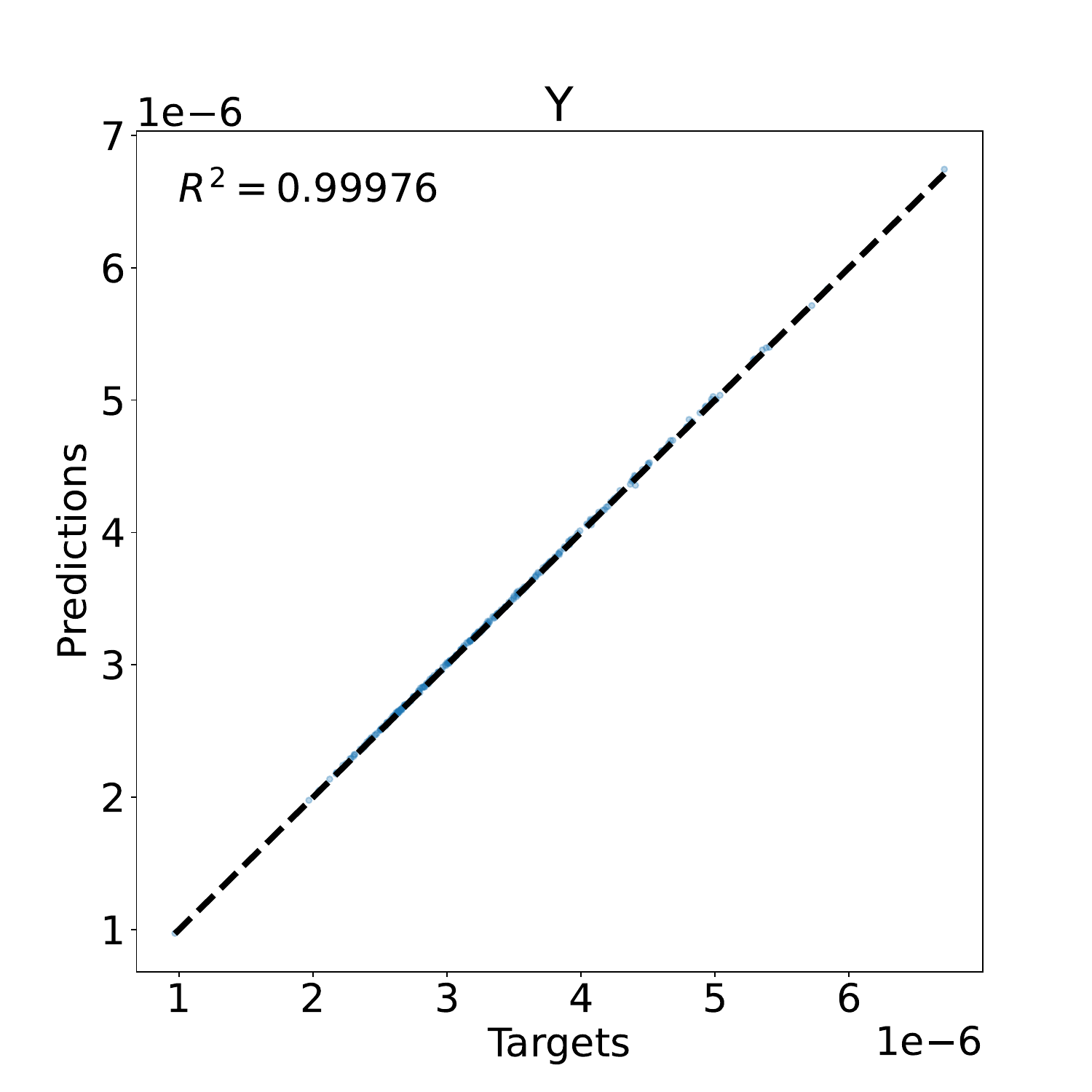}
    \caption{Surrogate A}
  \end{subfigure}
  \begin{subfigure}{0.4\textwidth}
    \includegraphics[width=\linewidth]{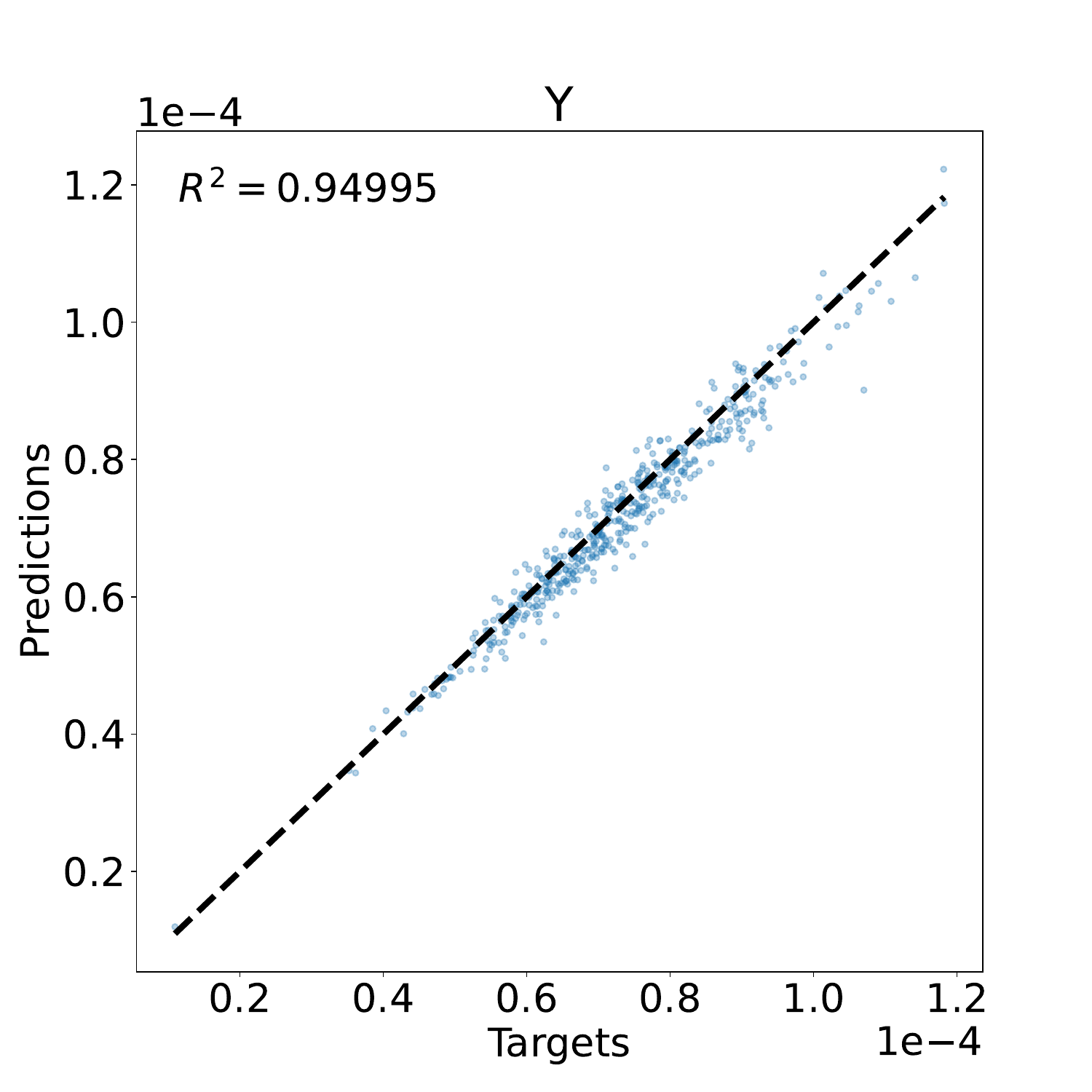}
    \caption{Surrogate B}
  \end{subfigure}
  \begin{subfigure}{0.4\textwidth}
    \includegraphics[width=\linewidth]{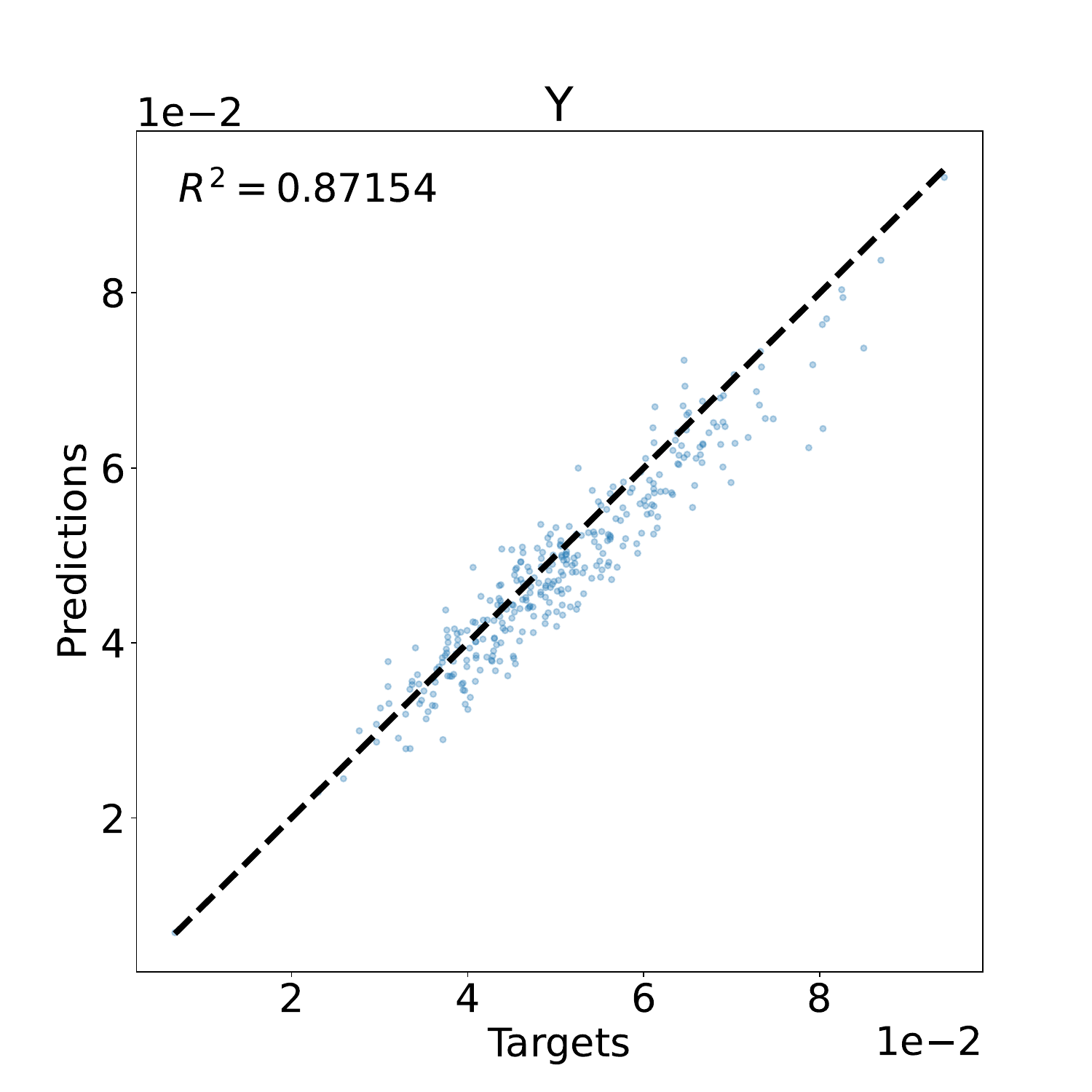}
    \caption{Surrogate C}
  \end{subfigure}
  \centering
  \caption{Comparison of outflow $Y$ from DFM models upscaled by numerical homogenization and by surrogates, $1000$~samples.}
  \label{fig:outflow_tar_pred}
\end{figure}

First, we compare upscaled DFM models on the Aquifer problem described in Section \ref{aquifer_problem}.
Fig. \ref{fig:outflow_tar_pred} presents the targets-predictions plot and $R^2$ for variable $Y$. For $K_f/K_m = \num{1e3}$, Surrogate A's predictions lead to almost the same $Y$ as using the numerical homogenization. In the case of $K_f/K_m = \num{1e5}$, Surrogate B's lower prediction accuracy causes a slightly lower $R^2 \approx 0.95$ for $Y$. Naturally, for $K_f/K_m = \num{1e7}$, the deteriorated prediction accuracy of Surrogate C leads to a somewhat smaller $R^2$ of about $0.87$ for $Y$ as well.

\begin{figure*}
  \centering
  \begin{subfigure}[t]{0.45\textwidth}
    \includegraphics[width=\linewidth]{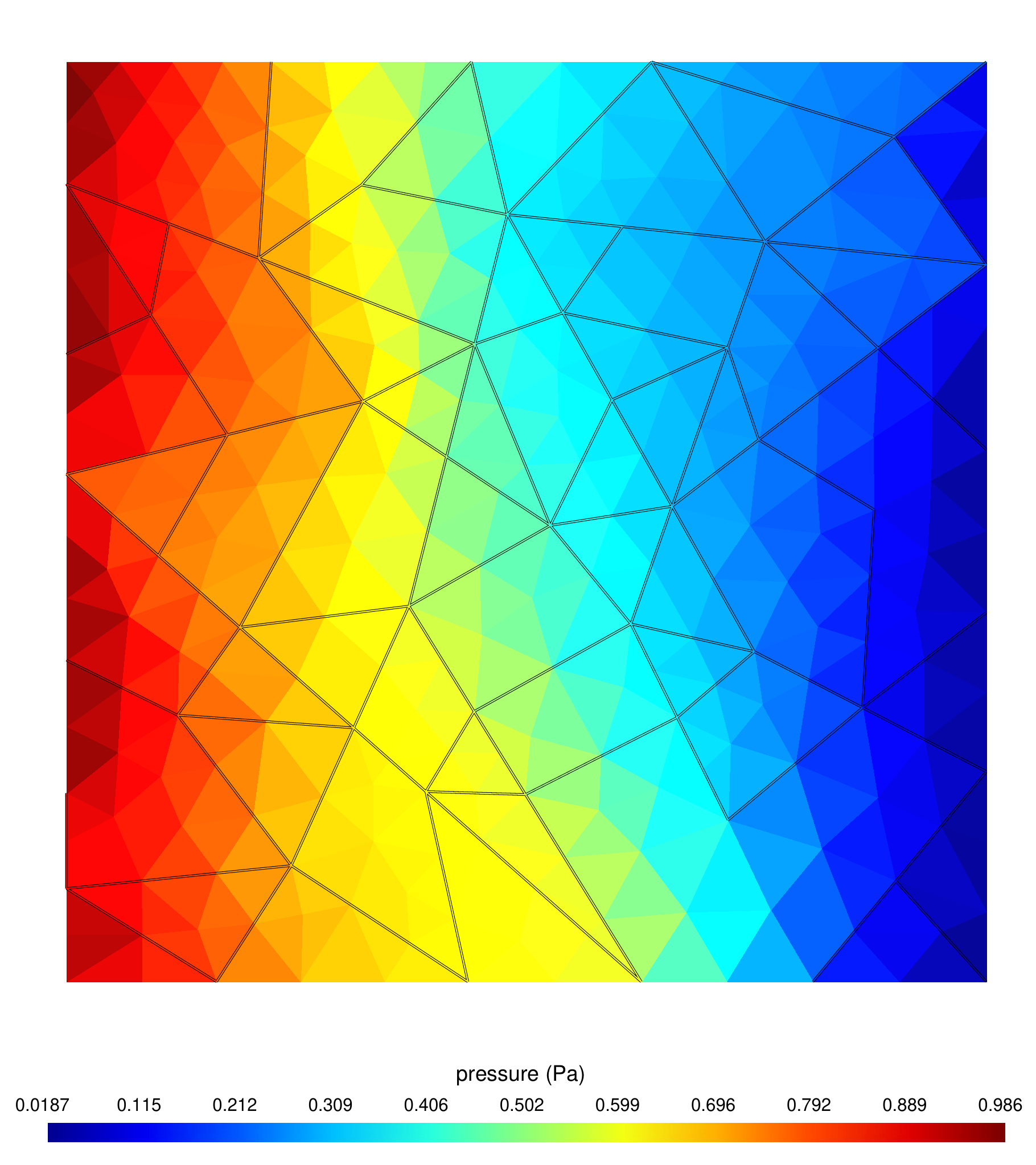}
    \caption{Upscaled DFM model by numerical homogenization}
  \end{subfigure}
  \hfill
  \begin{subfigure}[t]{0.45\textwidth}
    \includegraphics[width=\linewidth]{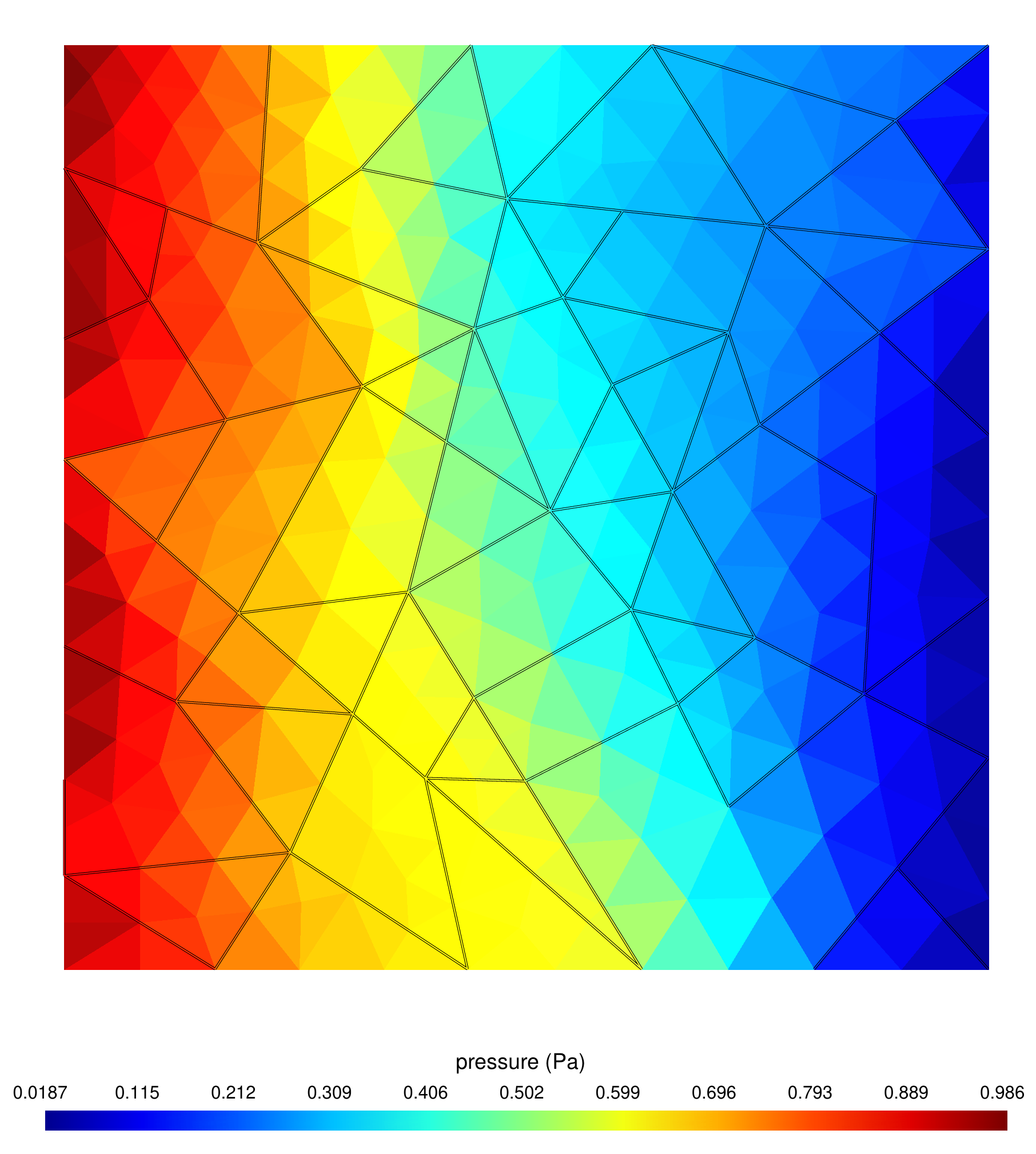}
    \caption{Upscaled DFM model by Surrogate A}
  \end{subfigure}
  \centering
  \caption{Aquifer problem - pressure field example, $K_f/K_m = \num{1e3}$}
  \label{fig:outflow_1e3}
\end{figure*}

\begin{figure*}
  \centering
  \begin{subfigure}[t]{0.45\textwidth}
    \includegraphics[width=\linewidth]{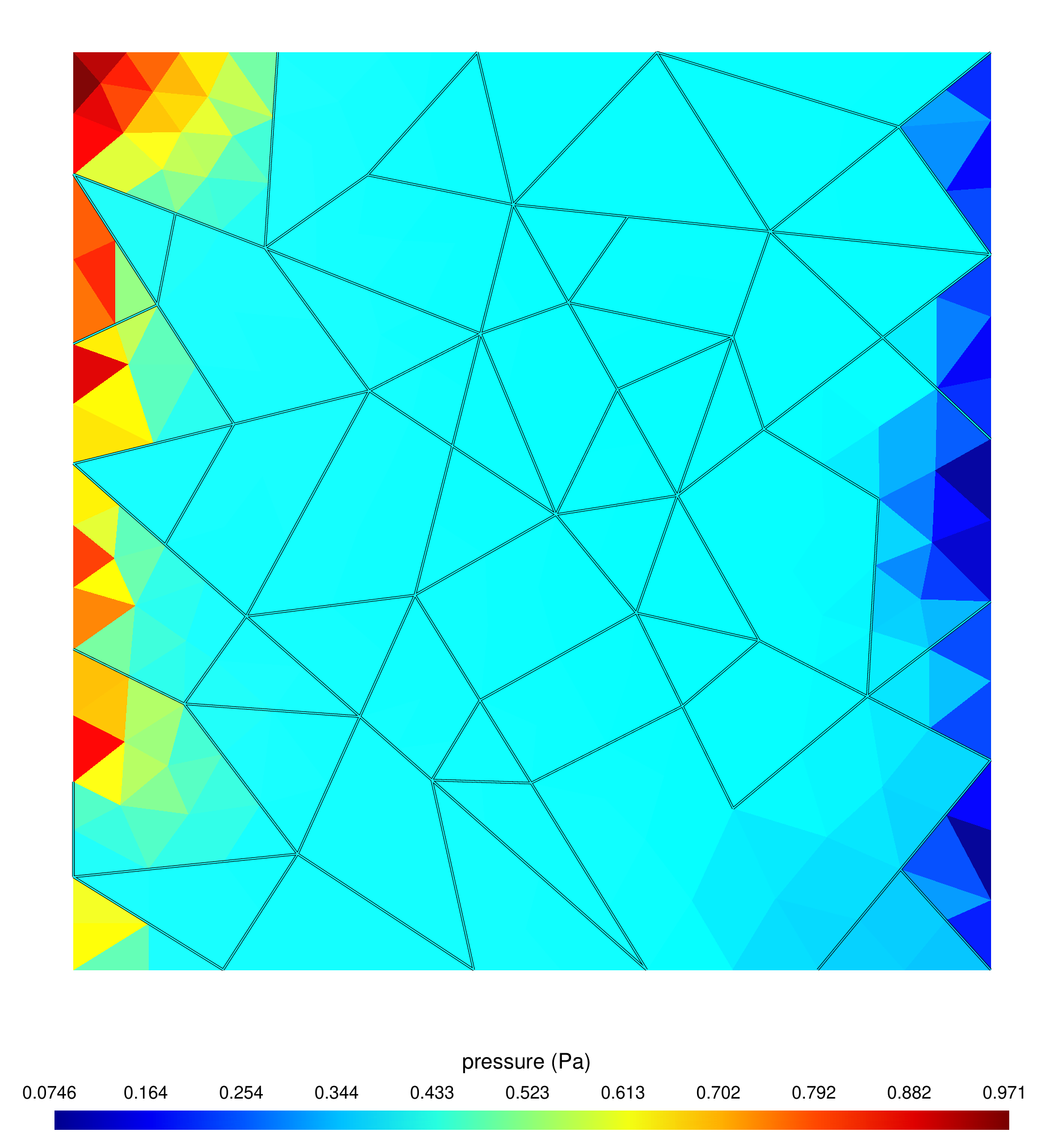}
    \caption{Upscaled DFM model by numerical homogenization}
  \end{subfigure}
  \hfill
  \begin{subfigure}[t]{0.445\textwidth}
    \includegraphics[width=\linewidth]{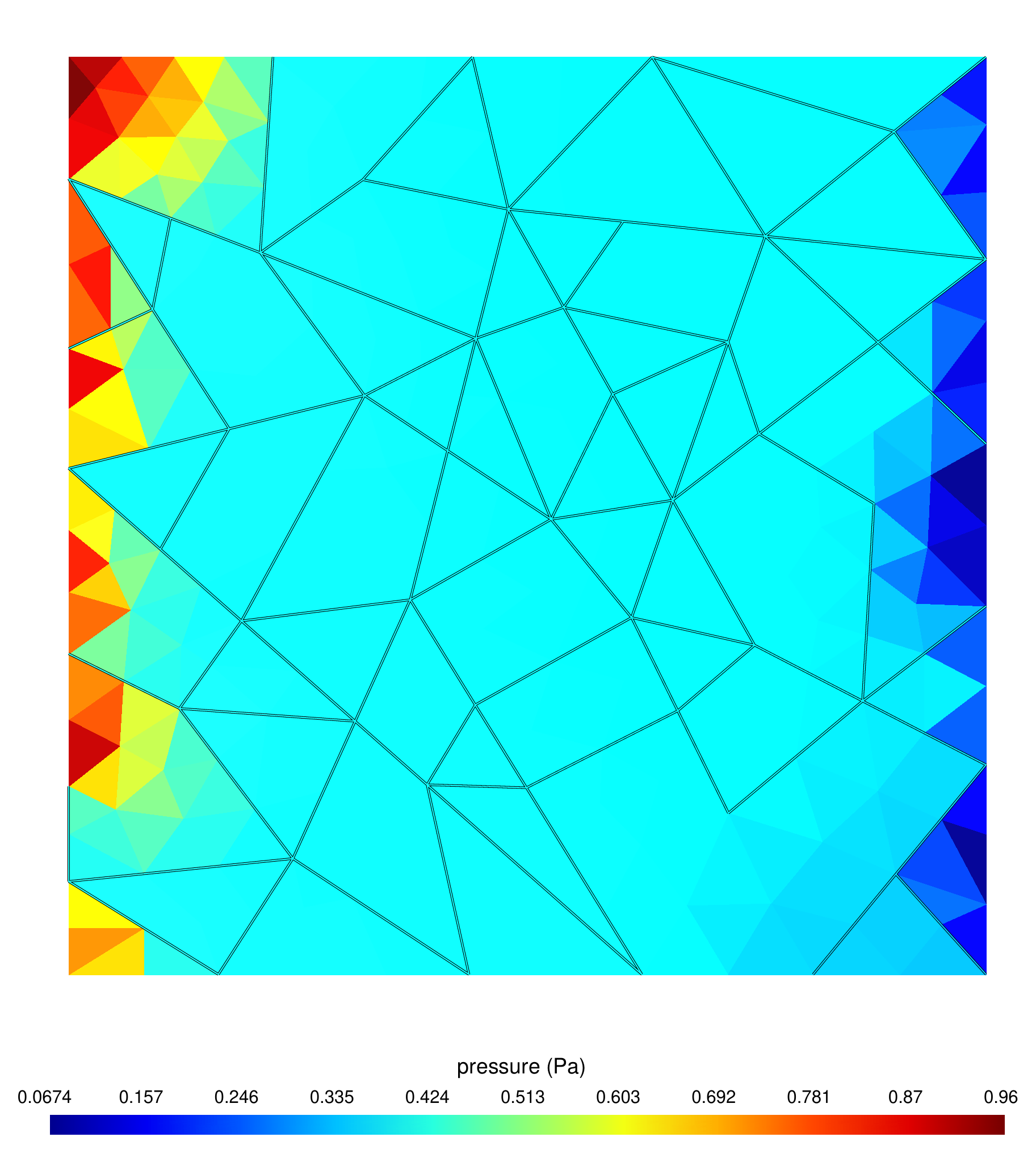}
    \caption{Upscaled DFM model by Surrogate C}
  \end{subfigure}
  \centering
  \caption{Aquifer problem - pressure field example, $K_f/K_m = \num{1e7}$}
  \label{fig:outflow_1e7}
\end{figure*}

While in the case of $K_f/K_m = \num{1e3}$, the whole upscaled SRF plays an equally important role in calculating $Y$ (see pressure fields in Fig. \ref{fig:outflow_1e3}), in the case of $K_f/K_m = \num{1e7}$ the influence of fractures is so dominant that mainly upscaled hydraulic conductivities near the left and right vertical boundaries are significant (see pressure fields in Fig. \ref{fig:outflow_1e7}). Therefore, the latter case is more sensitive to $\tn K^{eq}$ prediction accuracy on considerably fewer mesh elements.

\subsubsection{Anisotropy problem}
Second, we compare equivalent hydraulic conductivity tensors for DFM models upscaled by numerical homogenization and by our surrogates. For comparison, we calculate $\tn K^{eq}$ as described in Section \ref{anisotropy_problem} but for the entire domain $\Omega_0$ of the upscaled DFM model.
The obtained $R^2$ (see Table \ref{tab:equivalent_tensor_accuracy}) shows that all presented surrogates provide very accurate predictions of $\tn K^{eq}$, regardless of the value of $K_f/K_m \in \{\num{1e3}, \num{1e5}, \num{1e7}\}$ used. In other words, the components of $\tn K^{eq}$ of upscaled DFM models are nearly identical for numerical homogenization and surrogate modeling.

\begin{table}
    \centering
       \caption{$R^2$ comparison of components of calculated equivalent hydraulic conductivity tensors $\tn K^{eq}$ of DFM models upscaled by numerical homogenization and by surrogate}
    \begin{tabular}{lcccc}
        \toprule
        $R^2$ & $k_{xx}$ & $k_{xy}$ & $k_{yy}$  \\
        \midrule
        Surrogate A & $0.99967$ & $0.99990$ & $0.99910$ \\
        Surrogate B & $0.99960$ & $0.99984$ & $0.99989$ \\
        Surrogate C & $0.99998$ & $0.99998$ & $0.99998$  \\
        \bottomrule
    \end{tabular}
       \label{tab:equivalent_tensor_accuracy}
\end{table}

\begin{figure*}
  \centering
  \begin{subfigure}[t]{0.45\textwidth}
    \includegraphics[width=\linewidth]{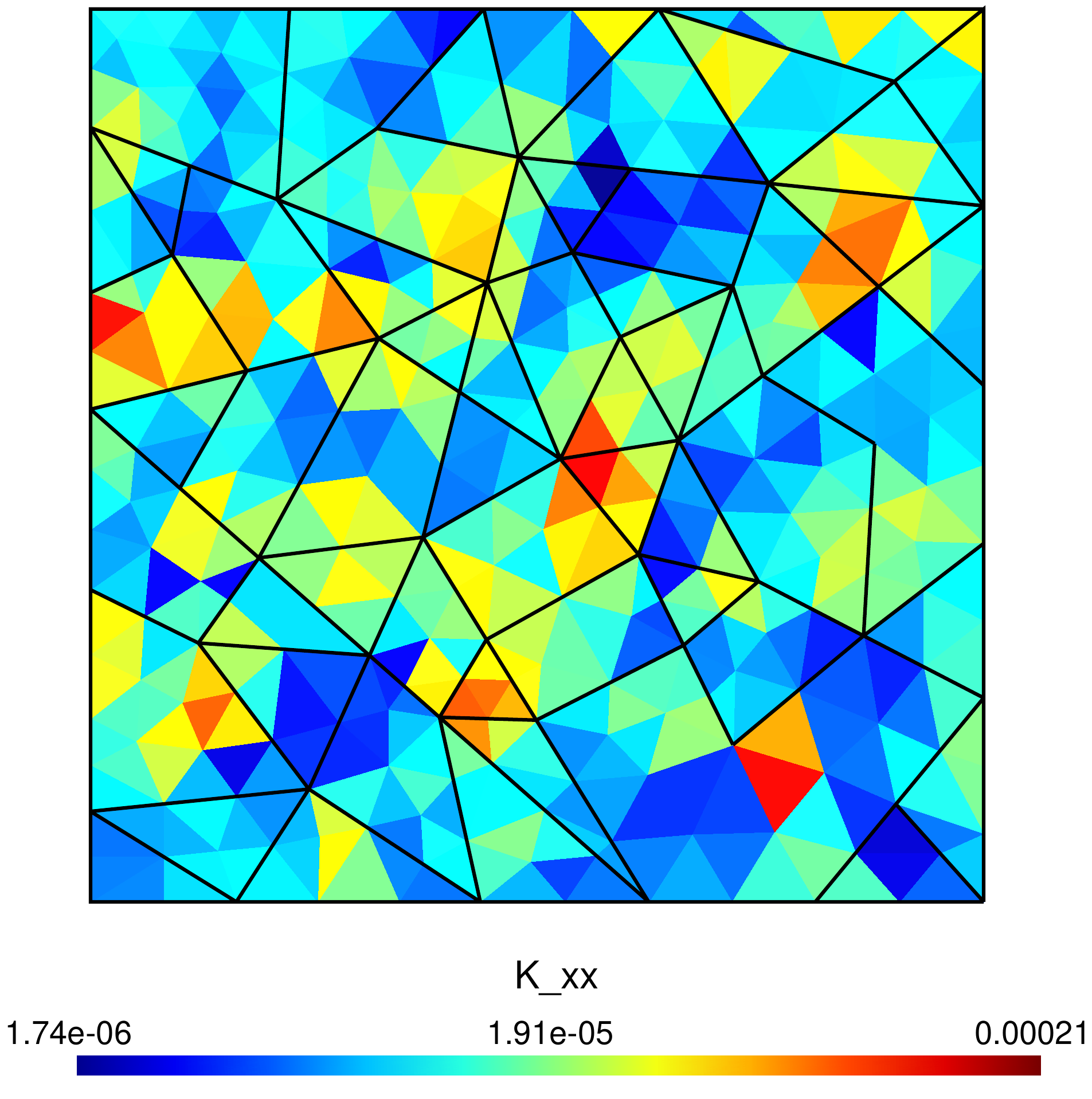}
    \caption{Upscaled SRF of $k_{xx}$ by numerical homogenization}
  \end{subfigure}
  \hfill
  \begin{subfigure}[t]{0.45\textwidth}
    \includegraphics[width=\linewidth]{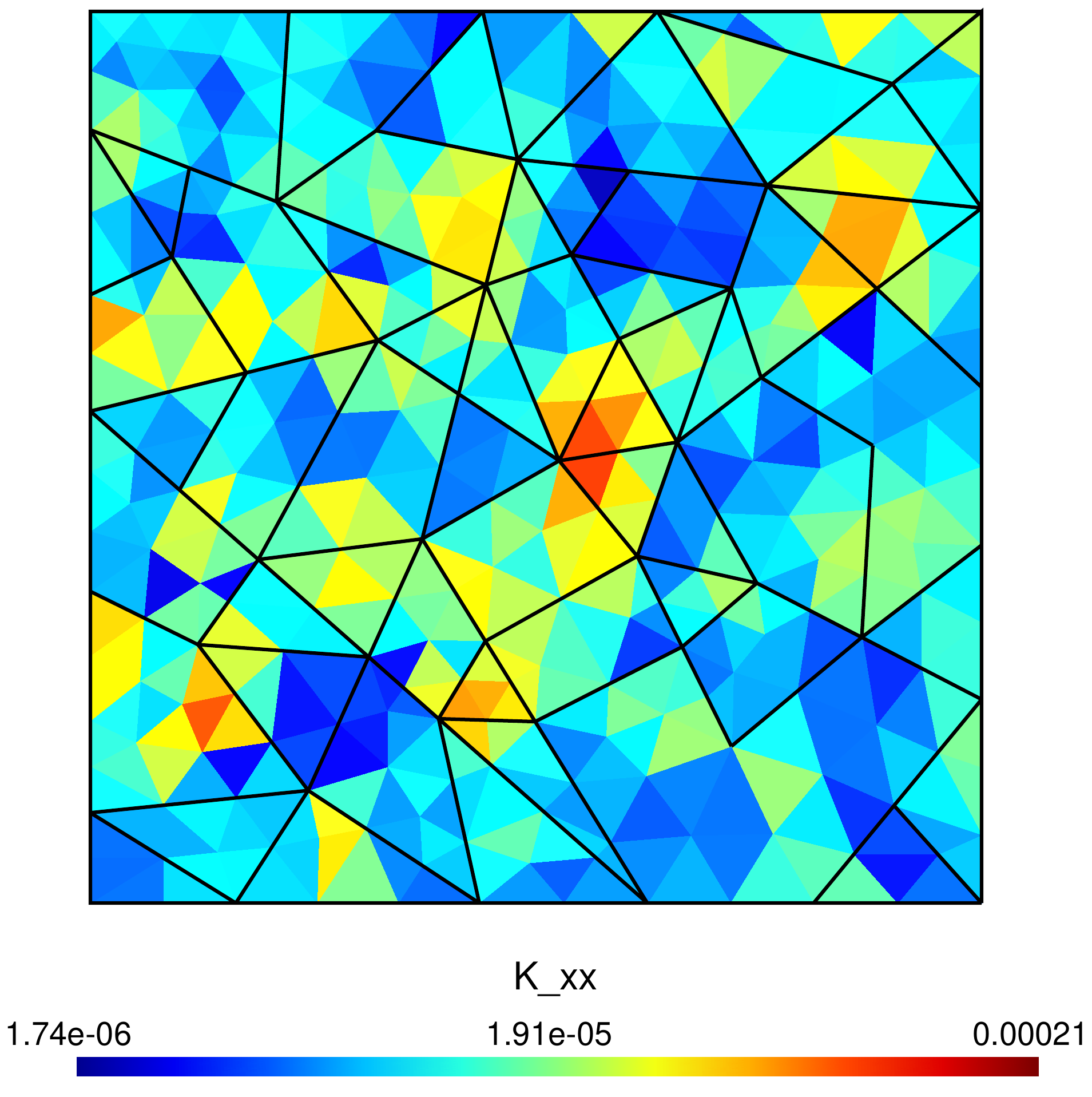}
    \caption{Upscaled SRF of $k_{xx}$ by Surrogate B}
  \end{subfigure}
  \centering
  \caption{Upscaled SRF of $k_{xx}$ by numerical homogenization and by Surrogate B, $K_f/K_m = \num{1e5}$}
  \label{fig:eff_tensor_fr_div_0_455_1_42}
\end{figure*}

In Fig. \ref{fig:eff_tensor_fr_div_0_455_1_42}, SRFs of $k_{xx}$ for an upscaled DFM model with $K_f/K_m = \num{1e5}$ are depicted. 
For the sake of clarity, hydraulic conductivity (K\_xx) scales are adjusted to have the same range for upscaling by Surrogate B and upscaling by numerical homogenization. The hydraulic conductivities of fractures are left out. 
Despite some dissimilarities, especially for higher values of $k_{xx}$, $\tn K^{eq}$ tensors are extremely similar. 
The explanation is that as the predictive accuracy of our surrogates decreases with increasing $K_f/K_m$, the impact of fracture conductivities increases simultaneously. Overall, we have highly accurate Surrogate A for $K_f/K_m = \num{1e3}$. Thus, $\tn K^{eq}$ is also highly accurate. For $K_f/K_m \in \{\num{1e5}, \num{1e7}\}$, the inaccuracy in predictions is compensated by the increased influence of fractures.

\section{Conclusions}\label{conclussion_section}
This article presented the deep learning-based surrogate for predicting the equivalent hydraulic conductivity tensor $\tn K^{eq}$ resulting from the upscaling of 2D DFM models. 
In addition to assessing the prediction accuracy of surrogates, we discussed their limitations and demonstrated practical usage on two macroscale problems.

We investigated different fracture-to-matrix hydraulic conductivity $K_f/K_m$ ratios and formed three separate datasets of the same structure for $K_f/K_m \in \{\num{1e3}, \num{1e5}, \num{1e7}\}$. For each ratio, a separate surrogate of the same architecture was trained and evaluated on test datasets.
Surrogates achieved highly accurate ($R^2 > 0.95$) predictions for all three independent components of $\tn K^{eq}$ when $K_f/K_m \leq 10^3$.
However, prediction accuracy deteriorated with increasing $K_f/K_m$ due to the complex distribution of $\tn K^{eq}$ components.
The off-diagonal component of $\tn K^{eq}$ exhibited consistently lower accuracy than diagonal components across all the cases. This inaccuracy is mainly driven by higher errors of predictions for data in the tails of the distributions. 
The speedup gained by surrogates depended on the number of homogenization blocks. From $4\times$ for $25$ blocks to $28\times$ for $1396$ blocks used for the upscaling of the smallest-scale DFM model in the study.

We also compared prediction accuracy for different fracture densities $\rho^{\prime}_{2D}$ and correlation lengths $\lambda $ of matrix SRF. In accordance with the literature, we observed improved accuracy on lower $\rho^{\prime}_{2D}$ than was explicitly used for training. Regarding correlation lengths, we noticed an improvement in prediction accuracy with increasing $\lambda$ for the case of lower impact of fractures ($K_f/K_m = \num{1e3}$) and conversely, a deterioration in prediction accuracy with increasing $\lambda$ for cases with the higher impact of fractures ($K_f/K_m \in \{\num{1e5}, \num{1e7}\}$). In all the cases, our surrogates exhibited expectedly accurate predictions even for $\lambda$ that was not included in the training set, up to $\lambda = 50$. In this matter, our surrogates generalize well.

Although the higher the prediction accuracy, the better, the necessary level of accuracy is not always obvious. Therefore, we investigated the impact of predicted  $\tn K^{eq}$ on two macroscale problems. 
In the case of comparing the usage of surrogates to numerical homogenization for upscaling of DFM models for the Aquifer problem, we observed $R^2 \gtrapprox 0.95$ for $K_f/K_m \in \{\num{1e3}, \num{1e5}\}$.  
Even the observed $R^2 \approx 0.87$ for $K_f/K_m = \num{1e7}$ might be sufficient. Such extreme hydraulic conductivity ratios may appear on coarse levels of MLMC, where the requested accuracy of the model is smaller.
Regarding the Anisotropy problem, the prediction capability of all of the presented surrogates was excellent ($R^2 \gtrapprox 0.999$) for upscaled DFM models.

We are also aware of the limitations of our approach and recognize the opportunities for improvement.
In rare cases where $K_f/K_m > \num{1e7}$ and $\rho^{\prime}_{2D} \leq 5$, prediction accuracy could be enhanced. 
Possible strategies include augmenting the training dataset in regions of smaller peaks and between peaks, as well as preprocessing $k_{xy}$ differently to ensure accurate prediction across all scales of values. 
For other considered $K_f/K_m$ ratios ($< \num{1e7}$), our surrogates demonstrated sufficient performance sufficiently for use within a multiscale MLMC scheme. This will be investigated in future work along with extending our approach to 3D DFM models.

\backmatter





\bmhead{Acknowledgements}
The research was supported within the framework of EURAD, the European Joint Programme on Radioactive Waste Management (Grant Agreement No 847593). 
This work was (partly) supported by the Student Grant Scheme at the Technical University of Liberec through project nr. SGS-2023-3379.
Computational resources were provided by the e-INFRA CZ project (ID:90254),
supported by the Ministry of Education, Youth and Sports of the Czech Republic.






\noindent 
\textbf{Code availability}
The source codes are available for downloading at the link:
\url{https://github.com/martinspetlik/MLMC-DFM/}

\noindent 
\textbf{Author contribution}
Martin Špetlík: Writing - original draft, Software, Experimentation.
Jan Březina: Supervision, Writing - review \& editing, Software.
Eric Laloy: Supervision, Writing - review \& editing.

\bibliography{sn-bibliography}


\begin{thebibliography}{52}
\ifx \bisbn   \undefined \def \bisbn  #1{ISBN #1}\fi
\ifx \binits  \undefined \def \binits#1{#1}\fi
\ifx \bauthor  \undefined \def \bauthor#1{#1}\fi
\ifx \batitle  \undefined \def \batitle#1{#1}\fi
\ifx \bjtitle  \undefined \def \bjtitle#1{#1}\fi
\ifx \bvolume  \undefined \def \bvolume#1{\textbf{#1}}\fi
\ifx \byear  \undefined \def \byear#1{#1}\fi
\ifx \bissue  \undefined \def \bissue#1{#1}\fi
\ifx \bfpage  \undefined \def \bfpage#1{#1}\fi
\ifx \blpage  \undefined \def \blpage #1{#1}\fi
\ifx \burl  \undefined \def \burl#1{\textsf{#1}}\fi
\ifx \doiurl  \undefined \def \doiurl#1{\url{https://doi.org/#1}}\fi
\ifx \betal  \undefined \def \betal{\textit{et al.}}\fi
\ifx \binstitute  \undefined \def \binstitute#1{#1}\fi
\ifx \binstitutionaled  \undefined \def \binstitutionaled#1{#1}\fi
\ifx \bctitle  \undefined \def \bctitle#1{#1}\fi
\ifx \beditor  \undefined \def \beditor#1{#1}\fi
\ifx \bpublisher  \undefined \def \bpublisher#1{#1}\fi
\ifx \bbtitle  \undefined \def \bbtitle#1{#1}\fi
\ifx \bedition  \undefined \def \bedition#1{#1}\fi
\ifx \bseriesno  \undefined \def \bseriesno#1{#1}\fi
\ifx \blocation  \undefined \def \blocation#1{#1}\fi
\ifx \bsertitle  \undefined \def \bsertitle#1{#1}\fi
\ifx \bsnm \undefined \def \bsnm#1{#1}\fi
\ifx \bsuffix \undefined \def \bsuffix#1{#1}\fi
\ifx \bparticle \undefined \def \bparticle#1{#1}\fi
\ifx \barticle \undefined \def \barticle#1{#1}\fi
\bibcommenthead
\ifx \bconfdate \undefined \def \bconfdate #1{#1}\fi
\ifx \botherref \undefined \def \botherref #1{#1}\fi
\ifx \url \undefined \def \url#1{\textsf{#1}}\fi
\ifx \bchapter \undefined \def \bchapter#1{#1}\fi
\ifx \bbook \undefined \def \bbook#1{#1}\fi
\ifx \bcomment \undefined \def \bcomment#1{#1}\fi
\ifx \oauthor \undefined \def \oauthor#1{#1}\fi
\ifx \citeauthoryear \undefined \def \citeauthoryear#1{#1}\fi
\ifx \endbibitem  \undefined \def \endbibitem {}\fi
\ifx \bconflocation  \undefined \def \bconflocation#1{#1}\fi
\ifx \arxivurl  \undefined \def \arxivurl#1{\textsf{#1}}\fi
\csname PreBibitemsHook\endcsname

\bibitem[\protect\citeauthoryear{Banks and Robins}{2002}]{Banks2002}
\begin{bbook}
\bauthor{\bsnm{Banks}, \binits{D.}},
\bauthor{\bsnm{Robins}, \binits{N.}}:
\bbtitle{An Introduction to Groundwater in Crystalline Bedrock}.
\bpublisher{Geological Survey of Norway},
\blocation{Trondheim, Norway}
(\byear{2002})
\end{bbook}
\endbibitem

\bibitem[\protect\citeauthoryear{Giles}{2015}]{Giles2015}
\begin{barticle}
\bauthor{\bsnm{Giles}, \binits{M.B.}}:
\batitle{Multilevel {{{M}onte {C}arlo}} methods}.
\bjtitle{Acta Numerica}
\bvolume{24},
\bfpage{259}--\blpage{328}
(\byear{2015})
\doiurl{10.1017/S096249291500001X}
\end{barticle}
\endbibitem

\bibitem[\protect\citeauthoryear{Auriault et~al.}{2009}]{Auriault20090101}
\begin{bbook}
\bauthor{\bsnm{Auriault}, \binits{J.-L.}},
\bauthor{\bsnm{Boutin}, \binits{C.}},
\bauthor{\bsnm{Geindreau}, \binits{C.}}:
\bbtitle{Homogenization of Coupled Phenomena in Heterogenous Media}.
\bpublisher{ISTE},
\blocation{London, UK}
(\byear{2009}).
\doiurl{10.1002/9780470612033}
\end{bbook}
\endbibitem

\bibitem[\protect\citeauthoryear{Bonnet et~al.}{2001}]{Bonnet2001ScalingOF}
\begin{barticle}
\bauthor{\bsnm{Bonnet}, \binits{E.}},
\bauthor{\bsnm{Bour}, \binits{O.}},
\bauthor{\bsnm{Odling}, \binits{N.E.}},
\bauthor{\bsnm{Davy}, \binits{P.}},
\bauthor{\bsnm{Main}, \binits{I.}},
\bauthor{\bsnm{Cowie}, \binits{P.}},
\bauthor{\bsnm{Berkowitz}, \binits{B.}}:
\batitle{Scaling of fracture systems in geological media}.
\bjtitle{Reviews of Geophysics}
\bvolume{39}(\bissue{3}),
\bfpage{347}--\blpage{383}
(\byear{2001})
\doiurl{10.1029/1999RG000074}
\end{barticle}
\endbibitem

\bibitem[\protect\citeauthoryear{Sanchez-Vila
  et~al.}{2006}]{https://doi.org/10.1029/2005RG000169}
\begin{botherref}
\oauthor{\bsnm{Sanchez-Vila}, \binits{X.}},
\oauthor{\bsnm{Guadagnini}, \binits{A.}},
\oauthor{\bsnm{Carrera}, \binits{J.}}:
Representative hydraulic conductivities in saturated groundwater flow.
Reviews of Geophysics
\textbf{44}(3)
(2006)
\doiurl{10.1029/2005RG000169}
\end{botherref}
\endbibitem

\bibitem[\protect\citeauthoryear{Renard and {de Marsily}}{1997}]{RENARD1997253}
\begin{barticle}
\bauthor{\bsnm{Renard}, \binits{P.}},
\bauthor{\bsnm{{de Marsily}}, \binits{G.}}:
\batitle{Calculating equivalent permeability: a review}.
\bjtitle{Advances in Water Resources}
\bvolume{20}(\bissue{5}),
\bfpage{253}--\blpage{278}
(\byear{1997})
\doiurl{10.1016/S0309-1708(96)00050-4}
\end{barticle}
\endbibitem

\bibitem[\protect\citeauthoryear{Chen et~al.}{2015}]{CHEN201560}
\begin{barticle}
\bauthor{\bsnm{Chen}, \binits{T.}},
\bauthor{\bsnm{Clauser}, \binits{C.}},
\bauthor{\bsnm{Marquart}, \binits{G.}},
\bauthor{\bsnm{Willbrand}, \binits{K.}},
\bauthor{\bsnm{Mottaghy}, \binits{D.}}:
\batitle{A new upscaling method for fractured porous media}.
\bjtitle{Advances in Water Resources}
\bvolume{80},
\bfpage{60}--\blpage{68}
(\byear{2015})
\doiurl{10.1016/j.advwatres.2015.03.009}
\end{barticle}
\endbibitem

\bibitem[\protect\citeauthoryear{Farmer}{2002}]{https://doi.org/10.1002/fld.267}
\begin{barticle}
\bauthor{\bsnm{Farmer}, \binits{C.L.}}:
\batitle{Upscaling: a review}.
\bjtitle{International Journal for Numerical Methods in Fluids}
\bvolume{40}(\bissue{1-2}),
\bfpage{63}--\blpage{78}
(\byear{2002})
\doiurl{10.1002/fld.267}
\end{barticle}
\endbibitem

\bibitem[\protect\citeauthoryear{Renard and Ababou}{2022}]{geosciences12070269}
\begin{botherref}
\oauthor{\bsnm{Renard}, \binits{P.}},
\oauthor{\bsnm{Ababou}, \binits{R.}}:
Equivalent permeability tensor of heterogeneous media: Upscaling methods and
  criteria (review and analyses).
Geosciences
\textbf{12}(7)
(2022)
\doiurl{10.3390/geosciences12070269}
\end{botherref}
\endbibitem

\bibitem[\protect\citeauthoryear{Bogdanov
  et~al.}{2003}]{https://doi.org/10.1029/2001WR000756}
\begin{botherref}
\oauthor{\bsnm{Bogdanov}, \binits{I.I.}},
\oauthor{\bsnm{Mourzenko}, \binits{V.V.}},
\oauthor{\bsnm{Thovert}, \binits{J.-F.}},
\oauthor{\bsnm{Adler}, \binits{P.M.}}:
Effective permeability of fractured porous media in steady state flow.
Water Resources Research
\textbf{39}(1)
(2003)
\doiurl{10.1029/2001WR000756}
\end{botherref}
\endbibitem

\bibitem[\protect\citeauthoryear{Koudina et~al.}{1998}]{Koudina1998}
\begin{barticle}
\bauthor{\bsnm{Koudina}, \binits{N.}},
\bauthor{\bsnm{Garcia}, \binits{R.G.}},
\bauthor{\bsnm{Thovert}, \binits{J.-F.}},
\bauthor{\bsnm{Adler}, \binits{P.M.}}:
\batitle{Permeability of three-dimensional fracture networks}.
\bjtitle{Physical Review E}
\bvolume{57}(\bissue{4}),
\bfpage{4466}--\blpage{4479}
(\byear{1998})
\doiurl{10.1103/PhysRevE.57.4466}
\end{barticle}
\endbibitem

\bibitem[\protect\citeauthoryear{Bogdanov et~al.}{2007}]{PhysRevE.76.036309}
\begin{barticle}
\bauthor{\bsnm{Bogdanov}, \binits{I.I.}},
\bauthor{\bsnm{Mourzenko}, \binits{V.V.}},
\bauthor{\bsnm{Thovert}, \binits{J.-F.}},
\bauthor{\bsnm{Adler}, \binits{P.M.}}:
\batitle{Effective permeability of fractured porous media with power-law
  distribution of fracture sizes}.
\bjtitle{Phys. Rev. E}
\bvolume{76},
\bfpage{036309}
(\byear{2007})
\doiurl{10.1103/PhysRevE.76.036309}
\end{barticle}
\endbibitem

\bibitem[\protect\citeauthoryear{Lang et~al.}{2014}]{Lang20140815}
\begin{barticle}
\bauthor{\bsnm{Lang}, \binits{P.S.}},
\bauthor{\bsnm{Paluszny}, \binits{A.}},
\bauthor{\bsnm{Zimmerman}, \binits{R.W.}}:
\batitle{Permeability tensor of three‐dimensional fractured porous rock and a
  comparison to trace map predictions}.
\bjtitle{Journal of Geophysical Research: Solid Earth}
\bvolume{119}(\bissue{8}),
\bfpage{6288}--\blpage{6307}
(\byear{2014})
\doiurl{10.1002/2014JB011027}
\end{barticle}
\endbibitem

\bibitem[\protect\citeauthoryear{Lee
  et~al.}{2001}]{https://doi.org/10.1029/2000WR900340}
\begin{barticle}
\bauthor{\bsnm{Lee}, \binits{S.H.}},
\bauthor{\bsnm{Lough}, \binits{M.F.}},
\bauthor{\bsnm{Jensen}, \binits{C.L.}}:
\batitle{Hierarchical modeling of flow in naturally fractured formations with
  multiple length scales}.
\bjtitle{Water Resources Research}
\bvolume{37}(\bissue{3}),
\bfpage{443}--\blpage{455}
(\byear{2001})
\doiurl{10.1029/2000WR900340}
\end{barticle}
\endbibitem

\bibitem[\protect\citeauthoryear{Azizmohammadi and
  Matthäi}{2017}]{Azizmohammadi2017}
\begin{barticle}
\bauthor{\bsnm{Azizmohammadi}, \binits{S.}},
\bauthor{\bsnm{Matthäi}, \binits{S.K.}}:
\batitle{Is the permeability of naturally fractured rocks scale dependent?}
\bjtitle{Water Resources Research}
\bvolume{53}(\bissue{9}),
\bfpage{8041}--\blpage{8063}
(\byear{2017})
\doiurl{10.1002/2016WR019764}
\end{barticle}
\endbibitem

\bibitem[\protect\citeauthoryear{Goodfellow
  et~al.}{2016}]{Goodfellow-et-al-2016}
\begin{bbook}
\bauthor{\bsnm{Goodfellow}, \binits{I.}},
\bauthor{\bsnm{Bengio}, \binits{Y.}},
\bauthor{\bsnm{Courville}, \binits{A.}}:
\bbtitle{Deep Learning}.
\bpublisher{MIT Press},
\blocation{Cambridge, Massachusetts}
(\byear{2016}).
\bcomment{\url{http://www.deeplearningbook.org}}
\end{bbook}
\endbibitem

\bibitem[\protect\citeauthoryear{Caglar et~al.}{2022}]{CAGLAR2022106973}
\begin{barticle}
\bauthor{\bsnm{Caglar}, \binits{B.}},
\bauthor{\bsnm{Broggi}, \binits{G.}},
\bauthor{\bsnm{Ali}, \binits{M.A.}},
\bauthor{\bsnm{Orgéas}, \binits{L.}},
\bauthor{\bsnm{Michaud}, \binits{V.}}:
\batitle{Deep learning accelerated prediction of the permeability of fibrous
  microstructures}.
\bjtitle{Composites Part A: Applied Science and Manufacturing}
\bvolume{158},
\bfpage{106973}
(\byear{2022})
\doiurl{10.1016/j.compositesa.2022.106973}
\end{barticle}
\endbibitem

\bibitem[\protect\citeauthoryear{Alqahtani et~al.}{2020}]{ALQAHTANI2020106514}
\begin{barticle}
\bauthor{\bsnm{Alqahtani}, \binits{N.}},
\bauthor{\bsnm{Alzubaidi}, \binits{F.}},
\bauthor{\bsnm{Armstrong}, \binits{R.T.}},
\bauthor{\bsnm{Swietojanski}, \binits{P.}},
\bauthor{\bsnm{Mostaghimi}, \binits{P.}}:
\batitle{Machine learning for predicting properties of porous media from 2d
  x-ray images}.
\bjtitle{Journal of Petroleum Science and Engineering}
\bvolume{184},
\bfpage{106514}
(\byear{2020})
\doiurl{10.1016/j.petrol.2019.106514}
\end{barticle}
\endbibitem

\bibitem[\protect\citeauthoryear{Wang et~al.}{2022}]{WANG2022123284}
\begin{barticle}
\bauthor{\bsnm{Wang}, \binits{Y.}},
\bauthor{\bsnm{Li}, \binits{H.}},
\bauthor{\bsnm{Xu}, \binits{J.}},
\bauthor{\bsnm{Liu}, \binits{S.}},
\bauthor{\bsnm{Wang}, \binits{X.}}:
\batitle{Machine learning assisted relative permeability upscaling for
  uncertainty quantification}.
\bjtitle{Energy}
\bvolume{245},
\bfpage{123284}
(\byear{2022})
\doiurl{10.1016/j.energy.2022.123284}
\end{barticle}
\endbibitem

\bibitem[\protect\citeauthoryear{Rao and Liu}{2020}]{RAO2020109850}
\begin{barticle}
\bauthor{\bsnm{Rao}, \binits{C.}},
\bauthor{\bsnm{Liu}, \binits{Y.}}:
\batitle{Three-dimensional convolutional neural network (3d-cnn) for
  heterogeneous material homogenization}.
\bjtitle{Computational Materials Science}
\bvolume{184},
\bfpage{109850}
(\byear{2020})
\doiurl{10.1016/j.commatsci.2020.109850}
\end{barticle}
\endbibitem

\bibitem[\protect\citeauthoryear{Hong and Liu}{2020}]{Hong2020}
\begin{barticle}
\bauthor{\bsnm{Hong}, \binits{J.}},
\bauthor{\bsnm{Liu}, \binits{J.}}:
\batitle{Rapid estimation of permeability from digital rock using 3d
  convolutional neural network}.
\bjtitle{Computational Geosciences}
\bvolume{24}(\bissue{4}),
\bfpage{1523}--\blpage{1539}
(\byear{2020})
\doiurl{10.1007/s10596-020-09941-w}
\end{barticle}
\endbibitem

\bibitem[\protect\citeauthoryear{Peng et~al.}{2022}]{Peng2022PhNetPM}
\begin{barticle}
\bauthor{\bsnm{Peng}, \binits{H.}},
\bauthor{\bsnm{Liu}, \binits{A.}},
\bauthor{\bsnm{Huang}, \binits{J.}},
\bauthor{\bsnm{Cao}, \binits{L.}},
\bauthor{\bsnm{Liu}, \binits{J.}},
\bauthor{\bsnm{Lu}, \binits{L.}}:
\batitle{Ph-net: Parallelepiped microstructure homogenization via 3d
  convolutional neural networks}.
\bjtitle{Additive Manufacturing}
\bvolume{60},
\bfpage{103237}
(\byear{2022})
\doiurl{10.1016/j.addma.2022.103237}
\end{barticle}
\endbibitem

\bibitem[\protect\citeauthoryear{Vasilyeva and
  Tyrylgin}{2021}]{VASILYEVA2021185}
\begin{barticle}
\bauthor{\bsnm{Vasilyeva}, \binits{M.}},
\bauthor{\bsnm{Tyrylgin}, \binits{A.}}:
\batitle{Machine learning for accelerating macroscopic parameters prediction
  for poroelasticity problem in stochastic media}.
\bjtitle{Computers \& Mathematics with Applications}
\bvolume{84},
\bfpage{185}--\blpage{202}
(\byear{2021})
\doiurl{10.1016/j.camwa.2020.09.024}
\end{barticle}
\endbibitem

\bibitem[\protect\citeauthoryear{Liu
  et~al.}{2023}]{https://doi.org/10.1029/2022JB025378}
\begin{barticle}
\bauthor{\bsnm{Liu}, \binits{M.}},
\bauthor{\bsnm{Ahmad}, \binits{R.}},
\bauthor{\bsnm{Cai}, \binits{W.}},
\bauthor{\bsnm{Mukerji}, \binits{T.}}:
\batitle{Hierarchical homogenization with deep-learning-based surrogate model
  for rapid estimation of effective permeability from digital rocks}.
\bjtitle{Journal of Geophysical Research: Solid Earth}
\bvolume{128}(\bissue{2}),
\bfpage{2022}--\blpage{025378}
(\byear{2023})
\doiurl{10.1029/2022JB025378}
\end{barticle}
\endbibitem

\bibitem[\protect\citeauthoryear{Meng et~al.}{2023}]{MENG2023104520}
\begin{barticle}
\bauthor{\bsnm{Meng}, \binits{Y.}},
\bauthor{\bsnm{Jiang}, \binits{J.}},
\bauthor{\bsnm{Wu}, \binits{J.}},
\bauthor{\bsnm{Wang}, \binits{D.}}:
\batitle{Transformer-based deep learning models for predicting permeability of
  porous media}.
\bjtitle{Advances in Water Resources}
\bvolume{179},
\bfpage{104520}
(\byear{2023})
\doiurl{10.1016/j.advwatres.2023.104520}
\end{barticle}
\endbibitem

\bibitem[\protect\citeauthoryear{Stepanov et~al.}{2023}]{STEPANOV2023114980}
\begin{barticle}
\bauthor{\bsnm{Stepanov}, \binits{S.}},
\bauthor{\bsnm{Spiridonov}, \binits{D.}},
\bauthor{\bsnm{Mai}, \binits{T.}}:
\batitle{Prediction of numerical homogenization using deep learning for the
  richards equation}.
\bjtitle{Journal of Computational and Applied Mathematics}
\bvolume{424},
\bfpage{114980}
(\byear{2023})
\doiurl{10.1016/j.cam.2022.114980}
\end{barticle}
\endbibitem

\bibitem[\protect\citeauthoryear{Pal et~al.}{2023}]{pr11020601}
\begin{botherref}
\oauthor{\bsnm{Pal}, \binits{M.}},
\oauthor{\bsnm{Makauskas}, \binits{P.}},
\oauthor{\bsnm{Malik}, \binits{S.}}:
Upscaling porous media using neural networks: A deep learning approach to
  homogenization and averaging.
Processes
\textbf{11}(2)
(2023)
\doiurl{10.3390/pr11020601}
\end{botherref}
\endbibitem

\bibitem[\protect\citeauthoryear{Cai et~al.}{2023}]{Cai_2023}
\begin{barticle}
\bauthor{\bsnm{Cai}, \binits{C.}},
\bauthor{\bsnm{Vlassis}, \binits{N.}},
\bauthor{\bsnm{Magee}, \binits{L.}},
\bauthor{\bsnm{Ma}, \binits{R.}},
\bauthor{\bsnm{Xiong}, \binits{Z.}},
\bauthor{\bsnm{Bahmani}, \binits{B.}},
\bauthor{\bsnm{Wong}, \binits{T.-F.}},
\bauthor{\bsnm{Wang}, \binits{Y.}},
\bauthor{\bsnm{Sun}, \binits{W.}}:
\batitle{Equivariant geometric learning for digital rock physics: estimating
  formation factor and effective permeability tensors from morse graph}.
\bjtitle{International Journal for Multiscale Computational Engineering}
\bvolume{21}(\bissue{5}),
\bfpage{1}--\blpage{24}
(\byear{2023})
\end{barticle}
\endbibitem

\bibitem[\protect\citeauthoryear{Ferreira et~al.}{2022}]{FERREIRA2022104264}
\begin{barticle}
\bauthor{\bsnm{Ferreira}, \binits{C.A.S.}},
\bauthor{\bsnm{Kadeethum}, \binits{T.}},
\bauthor{\bsnm{Bouklas}, \binits{N.}},
\bauthor{\bsnm{Nick}, \binits{H.M.}}:
\batitle{A framework for upscaling and modelling fluid flow for discrete
  fractures using conditional generative adversarial networks}.
\bjtitle{Advances in Water Resources}
\bvolume{166},
\bfpage{104264}
(\byear{2022})
\doiurl{10.1016/j.advwatres.2022.104264}
\end{barticle}
\endbibitem

\bibitem[\protect\citeauthoryear{He et~al.}{2021}]{10.2118/203901-MS}
\begin{bchapter}
\bauthor{\bsnm{He}, \binits{X.}},
\bauthor{\bsnm{Santoso}, \binits{R.}},
\bauthor{\bsnm{Alsinan}, \binits{M.}},
\bauthor{\bsnm{Kwak}, \binits{H.}},
\bauthor{\bsnm{Hoteit}, \binits{H.}}:
\bctitle{{Constructing Dual-Porosity Models from High-Resolution
  Discrete-Fracture Models Using Deep Neural Networks}}.
\bsertitle{SPE Reservoir Simulation Conference},
vol. \bseriesno{Day 1 Tue, October 26, 2021},
pp. \bfpage{011}--\blpage{014012}
(\byear{2021}).
\doiurl{10.2118/203901-MS}
\end{bchapter}
\endbibitem

\bibitem[\protect\citeauthoryear{Andrianov}{2022}]{Andrianov2022UpscalingOT}
\begin{barticle}
\bauthor{\bsnm{Andrianov}, \binits{N.}}:
\batitle{Upscaling of two-phase discrete fracture simulations using a
  convolutional neural network}.
\bjtitle{Computational Geosciences}
\bvolume{26},
\bfpage{1237}--\blpage{1259}
(\byear{2022})
\doiurl{10.1007/s10596-022-10149-3}
\end{barticle}
\endbibitem

\bibitem[\protect\citeauthoryear{Berre et~al.}{2019}]{Berre2019}
\begin{barticle}
\bauthor{\bsnm{Berre}, \binits{I.}},
\bauthor{\bsnm{Doster}, \binits{F.}},
\bauthor{\bsnm{Keilegavlen}, \binits{E.}}:
\batitle{Flow in fractured porous media}.
\bjtitle{Transport in Porous Media}
\bvolume{130}(\bissue{1}),
\bfpage{215}--\blpage{236}
(\byear{2019})
\doiurl{10.1007/s11242-018-1171-6}
\end{barticle}
\endbibitem

\bibitem[\protect\citeauthoryear{Long et~al.}{1982}]{Long1982Porousa}
\begin{barticle}
\bauthor{\bsnm{Long}, \binits{J.C.S.}},
\bauthor{\bsnm{Remer}, \binits{J.S.}},
\bauthor{\bsnm{Wilson}, \binits{C.R.}},
\bauthor{\bsnm{Witherspoon}, \binits{P.A.}}:
\batitle{Porous media equivalents for networks of discontinuous fractures}.
\bjtitle{Water Resources Research}
\bvolume{18}(\bissue{3}),
\bfpage{645}--\blpage{658}
(\byear{1982})
\doiurl{10.1029/WR018i003p00645}
\end{barticle}
\endbibitem

\bibitem[\protect\citeauthoryear{Bour and Davy}{1997}]{Bour1997Connectivity}
\begin{barticle}
\bauthor{\bsnm{Bour}, \binits{O.}},
\bauthor{\bsnm{Davy}, \binits{P.}}:
\batitle{Connectivity of random fault networks following a power law fault
  length distribution}.
\bjtitle{Water Resources Research}
\bvolume{33}(\bissue{7}),
\bfpage{1567}--\blpage{1583}
(\byear{1997})
\doiurl{10.1029/96WR00433}
\end{barticle}
\endbibitem

\bibitem[\protect\citeauthoryear{{de Dreuzy}
  et~al.}{2012}]{deDreuzy2012Influence}
\begin{botherref}
\oauthor{\bsnm{{de Dreuzy}}, \binits{J.-R.}},
\oauthor{\bsnm{M{\'e}heust}, \binits{Y.}},
\oauthor{\bsnm{Pichot}, \binits{G.}}:
Influence of fracture scale heterogeneity on the flow properties of
  three-dimensional discrete fracture networks ({{DFN}}).
Journal of Geophysical Research: Solid Earth
\textbf{117}(B11)
(2012)
\doiurl{10.1029/2012JB009461}
\end{botherref}
\endbibitem

\bibitem[\protect\citeauthoryear{Hardebol
  et~al.}{2015}]{https://doi.org/10.1002/2015JB011879}
\begin{barticle}
\bauthor{\bsnm{Hardebol}, \binits{N.J.}},
\bauthor{\bsnm{Maier}, \binits{C.}},
\bauthor{\bsnm{Nick}, \binits{H.}},
\bauthor{\bsnm{Geiger}, \binits{S.}},
\bauthor{\bsnm{Bertotti}, \binits{G.}},
\bauthor{\bsnm{Boro}, \binits{H.}}:
\batitle{Multiscale fracture network characterization and impact on flow: A
  case study on the latemar carbonate platform}.
\bjtitle{Journal of Geophysical Research: Solid Earth}
\bvolume{120}(\bissue{12}),
\bfpage{8197}--\blpage{8222}
(\byear{2015})
\doiurl{10.1002/2015JB011879}
\end{barticle}
\endbibitem

\bibitem[\protect\citeauthoryear{Adler and Thovert}{1999}]{Adler1999}
\begin{bbook}
\bauthor{\bsnm{Adler}, \binits{P.M.}},
\bauthor{\bsnm{Thovert}, \binits{J.-F.}}:
\bbtitle{Fractures and Fracture Networks}.
\bpublisher{Springer},
\blocation{Dordrecht}
(\byear{1999}).
\doiurl{10.1007/978-94-017-1599-7}
\end{bbook}
\endbibitem

\bibitem[\protect\citeauthoryear{Sahimi}{2011}]{Sahimi20110420}
\begin{bbook}
\bauthor{\bsnm{Sahimi}, \binits{M.}}:
\bbtitle{Flow and Transport in Porous Media and Fractured Rock},
\bedition{2.} edn.
\bpublisher{Wiley},
\blocation{Germany}
(\byear{2011}).
\doiurl{10.1002/9783527636693}
\end{bbook}
\endbibitem

\bibitem[\protect\citeauthoryear{Liu et~al.}{2016}]{Liu2016}
\begin{barticle}
\bauthor{\bsnm{Liu}, \binits{R.}},
\bauthor{\bsnm{Li}, \binits{B.}},
\bauthor{\bsnm{Jiang}, \binits{Y.}},
\bauthor{\bsnm{Huang}, \binits{N.}}:
\batitle{Review: Mathematical expressions for estimating equivalent
  permeability of rock fracture networks}.
\bjtitle{Hydrogeology Journal}
\bvolume{24}(\bissue{7}),
\bfpage{1623}--\blpage{1649}
(\byear{2016})
\doiurl{10.1007/s10040-016-1441-8}
\end{barticle}
\endbibitem

\bibitem[\protect\citeauthoryear{Hadgu et~al.}{2017}]{Hadgu_comparative_2017}
\begin{barticle}
\bauthor{\bsnm{Hadgu}, \binits{T.}},
\bauthor{\bsnm{Karra}, \binits{S.}},
\bauthor{\bsnm{Kalinina}, \binits{E.}},
\bauthor{\bsnm{Makedonska}, \binits{N.}},
\bauthor{\bsnm{Hyman}, \binits{J.D.}},
\bauthor{\bsnm{Klise}, \binits{K.}},
\bauthor{\bsnm{Viswanathan}, \binits{H.S.}},
\bauthor{\bsnm{Wang}, \binits{Y.}}:
\batitle{A comparative study of discrete fracture network and equivalent
  continuum models for simulating flow and transport in the far field of a
  hypothetical nuclear waste repository in crystalline host rock}.
\bjtitle{Journal of Hydrology}
\bvolume{553},
\bfpage{59}--\blpage{70}
(\byear{2017})
\doiurl{10.1016/j.jhydrol.2017.07.046}
\end{barticle}
\endbibitem

\bibitem[\protect\citeauthoryear{Kottwitz
  et~al.}{2021}]{Kottwitz_Investigating_2021}
\begin{barticle}
\bauthor{\bsnm{Kottwitz}, \binits{M.O.}},
\bauthor{\bsnm{Popov}, \binits{A.A.}},
\bauthor{\bsnm{Abe}, \binits{S.}},
\bauthor{\bsnm{Kaus}, \binits{B.J.P.}}:
\batitle{Investigating the effects of intersection flow localization in
  equivalent-continuum-based upscaling of flow in discrete fracture networks}.
\bjtitle{Solid Earth}
\bvolume{12}(\bissue{10}),
\bfpage{2235}--\blpage{2254}
(\byear{2021})
\doiurl{10.5194/se-12-2235-2021}
\end{barticle}
\endbibitem

\bibitem[\protect\citeauthoryear{Müller and Schüler}{2019}]{GSTools}
\begin{botherref}
\oauthor{\bsnm{Müller}, \binits{S.}},
\oauthor{\bsnm{Schüler}, \binits{L.}}:
{GSTools}.
\url{https://github.com/GeoStat-Framework/GSTools}
(2019)
\end{botherref}
\endbibitem

\bibitem[\protect\citeauthoryear{Sandve et~al.}{2012}]{Sandve_efficient_2012a}
\begin{barticle}
\bauthor{\bsnm{Sandve}, \binits{T.H.}},
\bauthor{\bsnm{Berre}, \binits{I.}},
\bauthor{\bsnm{Nordbotten}, \binits{J.M.}}:
\batitle{An efficient multi-point flux approximation method for {{Discrete
  Fracture}}--{{Matrix}} simulations}.
\bjtitle{Journal of Computational Physics}
\bvolume{231}(\bissue{9}),
\bfpage{3784}--\blpage{3800}
(\byear{2012})
\doiurl{10.1016/j.jcp.2012.01.023}
\end{barticle}
\endbibitem

\bibitem[\protect\citeauthoryear{Berrone
  et~al.}{2013}]{Berrone_Simulations_2013}
\begin{barticle}
\bauthor{\bsnm{Berrone}, \binits{S.}},
\bauthor{\bsnm{Pieraccini}, \binits{S.}},
\bauthor{\bsnm{Scial{\`o}}, \binits{S.}}:
\batitle{On {{Simulations}} of {{Discrete Fracture Network Flows}} with an
  {{Optimization-Based Extended Finite Element Method}}}.
\bjtitle{SIAM Journal on Scientific Computing}
\bvolume{35}(\bissue{2}),
\bfpage{908}--\blpage{935}
(\byear{2013})
\doiurl{10.1137/120882883}
\end{barticle}
\endbibitem

\bibitem[\protect\citeauthoryear{B{\v{r}}ezina and
  Stebel}{2016}]{Brezina2016Analysis}
\begin{bchapter}
\bauthor{\bsnm{B{\v{r}}ezina}, \binits{J.}},
\bauthor{\bsnm{Stebel}, \binits{J.}}:
\bctitle{Analysis of model error for a continuum-fracture model of porous media
  flow}.
In: \beditor{\bsnm{Kozubek}, \binits{T.}},
\beditor{\bsnm{Blaheta}, \binits{R.}},
\beditor{\bsnm{{\v{S}}{\'i}stek}, \binits{J.}},
\beditor{\bsnm{Rozlo{\v{z}}n{\'i}k}, \binits{M.}},
\beditor{\bsnm{{\v{C}}erm{\'a}k}, \binits{M.}} (eds.)
\bbtitle{High Performance Computing in Science and Engineering},
pp. \bfpage{152}--\blpage{160}.
\bpublisher{Springer},
\blocation{Cham}
(\byear{2016}).
\doiurl{10.1007/978-3-319-40361-8_11}
\end{bchapter}
\endbibitem

\bibitem[\protect\citeauthoryear{B{\v r}ezina et~al.}{2011--2023}]{flow123d}
\begin{botherref}
\oauthor{\bsnm{B{\v r}ezina}, \binits{J.}},
\oauthor{\bsnm{Stebel}, \binits{J.}},
\oauthor{\bsnm{Exner}, \binits{P.}},
\oauthor{\bsnm{Hyb{\v s}}, \binits{J.}}:
Flow123d.
\url{http://flow123d.github.com}, repository:
  \url{http://github.com/flow123d/flow123d}
(2011--2023)
\end{botherref}
\endbibitem

\bibitem[\protect\citeauthoryear{Durlofsky}{1991}]{Durlofsky1991}
\begin{barticle}
\bauthor{\bsnm{Durlofsky}, \binits{L.J.}}:
\batitle{Numerical calculation of equivalent grid block permeability tensors
  for heterogeneous porous media}.
\bjtitle{Water Resources Research}
\bvolume{27}(\bissue{5}),
\bfpage{699}--\blpage{708}
(\byear{1991})
\doiurl{10.1029/91WR00107}
\end{barticle}
\endbibitem

\bibitem[\protect\citeauthoryear{Alzubaidi et~al.}{2021}]{Alzubaidi2021}
\begin{botherref}
\oauthor{\bsnm{Alzubaidi}, \binits{L.}},
\oauthor{\bsnm{Zhang}, \binits{J.}},
\oauthor{\bsnm{Humaidi}, \binits{A.J.}},
\oauthor{\bsnm{Al-Dujaili}, \binits{A.}},
\oauthor{\bsnm{Duan}, \binits{Y.}},
\oauthor{\bsnm{Al-Shamma}, \binits{O.}},
\oauthor{\bsnm{Santamaría}, \binits{J.}},
\oauthor{\bsnm{Fadhel}, \binits{M.A.}},
\oauthor{\bsnm{Al-Amidie}, \binits{M.}},
\oauthor{\bsnm{Farhan}, \binits{L.}}:
Review of deep learning.
Journal of Big Data
\textbf{8}(1)
(2021)
\doiurl{10.1186/s40537-021-00444-8}
\end{botherref}
\endbibitem

\bibitem[\protect\citeauthoryear{Zhang et~al.}{2020}]{zhang2020dive}
\begin{bbook}
\bauthor{\bsnm{Zhang}, \binits{A.}},
\bauthor{\bsnm{Lipton}, \binits{Z.C.}},
\bauthor{\bsnm{Li}, \binits{M.}},
\bauthor{\bsnm{Smola}, \binits{A.J.}}:
\bbtitle{Dive Into Deep Learning},
(\byear{2020}).
\bcomment{\url{https://d2l.ai}}
\end{bbook}
\endbibitem

\bibitem[\protect\citeauthoryear{Bednar et~al.}{2022}]{datashader}
\begin{botherref}
\oauthor{\bsnm{Bednar}, \binits{J.A.}},
\oauthor{\bsnm{Crist}, \binits{J.}},
\oauthor{\bsnm{Cottam}, \binits{J.}},
\oauthor{\bsnm{Wang}, \binits{P.}}:
{Datashader}.
\url{https://datashader.org/}
(2022)
\end{botherref}
\endbibitem

\bibitem[\protect\citeauthoryear{He et~al.}{2023}]{HE2023592}
\begin{barticle}
\bauthor{\bsnm{He}, \binits{C.}},
\bauthor{\bsnm{Yao}, \binits{C.}},
\bauthor{\bsnm{Jin}, \binits{Y.-z.}},
\bauthor{\bsnm{Jiang}, \binits{Q.-h.}},
\bauthor{\bsnm{Zhou}, \binits{C.-b.}}:
\batitle{Effective permeability of fractured porous media with fracture density
  near the percolation threshold}.
\bjtitle{Applied Mathematical Modelling}
\bvolume{117},
\bfpage{592}--\blpage{608}
(\byear{2023})
\doiurl{10.1016/j.apm.2023.01.010}
\end{barticle}
\endbibitem

\bibitem[\protect\citeauthoryear{Zhu et~al.}{2023}]{ZHU2023111186}
\begin{barticle}
\bauthor{\bsnm{Zhu}, \binits{C.}},
\bauthor{\bsnm{Wang}, \binits{J.}},
\bauthor{\bsnm{Sang}, \binits{S.}},
\bauthor{\bsnm{Liang}, \binits{W.}}:
\batitle{A multiscale neural network model for the prediction on the equivalent
  permeability of discrete fracture network}.
\bjtitle{Journal of Petroleum Science and Engineering}
\bvolume{220},
\bfpage{111186}
(\byear{2023})
\doiurl{10.1016/j.petrol.2022.111186}
\end{barticle}
\endbibitem

\end{thebibliography}

\end{document}